\documentclass[10pt,a4paper]{amsart}
\usepackage[latin1]{inputenc}
\usepackage{amsmath,amssymb,amsthm}
\usepackage{graphicx}
\usepackage[left=2.5cm,right=2.5cm,top=2.5cm,bottom=2.5cm]{geometry}
\usepackage{mathtools}
\usepackage{color}
\usepackage{hyperref}

\def\C{\mathbb{C}}
\def\R{\mathbb{R}}
\def\Z{\mathbb{Z}}
\def\N{\mathbb{N}}

\def\dx{\Delta x}
\def\dt{\Delta t}

\def\ud{u_\Delta}
\def\umu{u^\mu}
\def\vmu{v^\mu}
\def\wmu{w^\mu}

\def\coefzero{F^\Delta_0}
\def\coefone{F^\Delta_1}
\def\coeftwo{F^\Delta_2}

\def\sign{\mathop{\mathrm{sign}}}
\newcommand{\squarehat}[1]{\overbracket[.7pt][.5ex]{#1}}

\newtheorem{thrm}{Theorem}[section]
\newtheorem{prpstn}{Proposition}[section]
\newtheorem{lmm}{Lemma}[section]
\newtheorem{rmrk}{Remark}

\begin{document}

\title[A semi-discrete scheme for the augmented Burgers equation]{A semi-discrete large-time behavior preserving scheme \\ for the augmented Burgers equation}
\author[L. I. Ignat, A. Pozo]{L. I. Ignat, A. Pozo}
\date{\footnotesize \today}

\address{Liviu I. Ignat
\hfill\break\indent Institute of Mathematics ``Simion Stoilow'' of the Romanian Academy\\
\hfill\break\indent 21 Calea Grivitei Street, 010702 Bucharest, Romania.
}
 \email{{\tt liviu.ignat@gmail.com}\hfill\break\indent {\it Web page: }{\tt http://www.imar.ro/\~{}lignat}}
\address{Alejandro Pozo
\hfill\break\indent Innovalia Association\\
\hfill\break\indent Carretera de Asua 6, 48930 Las Arenas - Getxo, Spain\\
\hfill\break\indent \and
\hfill\break\indent BCAM - Basque Center for Applied Mathematics\\
\hfill\break\indent Alameda de Mazarredo 14, 48009 Bilbao, Spain.}
 \email{{\tt alejandropozo@gmail.com}}
 

\date{June 1st, 2017}

\begin{abstract}
	In this paper we analyze the large-time behavior of the augmented Burgers equation. We first study the well-posedness of the Cauchy problem and obtain $L^1$-$L^p$ decay rates. The asymptotic behavior of the solution is obtained by showing that the influence of the convolution term $K*u_{xx}$ is the same  as $u_{xx}$ for large times. Then, we propose a semi-discrete numerical scheme that preserves this asymptotic behavior, by introducing two correcting factors in the discretization of the non-local term. Numerical experiments illustrating the accuracy of the results of the paper are also presented.
\end{abstract}

\maketitle

%%%%%%%%%%%% SECTION 1 %%%%%%%%%%%%%%%%%%%%%%%%%

\section{Introduction and main results}\label{section1}
In this paper we consider the following equation:
\begin{equation}\label{eq:modelfull}
	\begin{cases}
		u_t=uu_x+ \nu\, u_{xx}+ c\, K_\theta*u_{xx},&(t,x)\in(0,\infty)\times\R,\\
		u(0,x)=u_0(x),&x\in\R,
	\end{cases}
\end{equation}
where $*$ denotes the convolution in the $x$ variable, the parameters $\nu,c,\theta$ are positive and
\begin{equation}\label{eq:kernelk}
	K_\theta(z)=\begin{cases} \frac1\theta e^{-z/\theta},&z>0, \\ 0,& \mbox{elsewhere.} \end{cases}
\end{equation}
This is a constant-parameter version of the augmented Burgers equation, which has been used to model the propagation of the sonic-boom produced by supersonic aircrafts from their near-field down to the ground level.

Until the last decade of the 20th century, linear theory was used to model the evolution of this strident noise, based on the seminal works by Hayes \cite{Hayes:1947} and Whitham \cite{Whitham:1952}. In fact, until recently, most of  the research, both from and analytical and a numerical point of view, followed the so-called Jones-Seebass-George-Darden theory for sonic-boom minimization \cite{Darden:1977,Jones:1961,Seebass:1969a,Seebass:1969b,Seebass:1972}.

Newer trends have started to use nonlinear physical models to improve the characterization of the sonic-boom propagation. In this paper we focus on Burgers-type equations, which have been one of the main tools to model the propagation of finite-amplitude plane waves. The classical viscous Burgers equation \cite{Burgers:1940} was first considered for wave propagation in a lossy medium. Successive generalizations included other effects such as geometrical spreading and inhomogeneous mediums \cite{Carlton:1974,Fridman:1976,Lighthill:1956} or relaxation processes \cite{Coulouvrat:2009,Pierce:1989}. All those phenomena were taken into account in the augmented Burgers equation, initially developed by Cleveland \cite{Cleveland:1995} and then adopted by Rallabhandi \cite{Rallabhandi:2011a, Rallabhandi:2011b}. This equation is given by
\begin{equation}\label{eq:sonic}
	\frac{\partial P}{\partial \sigma} = P \frac{\partial P}{\partial \tau} + \frac{1}{\Gamma}\frac{\partial^2 P}{\partial \tau^2}+\sum_{\nu} C_\nu\frac{1}{1+\theta_\nu\frac{\partial}{\partial \tau}}\frac{\partial^2P}{\partial\tau^2}-\frac{1}{2G}\frac{\partial G}{\partial\sigma}P+\frac{1}{2\rho_0c_0}\frac{\partial(\rho_0c_0)}{\partial\sigma} P,
\end{equation}
where $P=P(\sigma,\tau)$ is the dimensionless perturbation of the pressure distribution. The covered distance $\sigma$ and time of the perturbation $\tau$ are also dimensionless. The operator appearing in the summation, corresponding to the molecular relaxations, it is defined by:
\begin{equation}\label{eq:operator}
	\frac{1}{1+\theta_\nu\frac{\partial}{\partial\tau}} f(\tau) = \frac{1}{\theta_\nu}\int_{-\infty}^\tau e^{(\xi-\tau)/\theta_\nu} f(\xi) d\xi = K_{\theta_\nu}*f(\tau),
\end{equation}
Typically, two relaxation modes are considered: one for Oxygen molecules and another one for Nitrogen ones. $\theta_\nu$ and $C_\nu$ are the dimensionless relaxation time and dispersion parameter, respectively, for each one. $\Gamma$ is a dimensionless thermo-viscous parameter and function $G\equiv G(\sigma)$ denotes the ray-tube area. The atmosphere conditions are given by density $\rho_0\equiv\rho(\sigma)$ and speed of sound $c_0\equiv c_0(\sigma)$, both closely related to the altitude of the flight. We refer the reader to \cite{Cleveland:1995} for a detailed description on the development of this model and to \cite{Alonso:2012} for a comprehensive review about the sonic-boom minimization problem.

Industrial applications of this kind of models, such as the aforementioned sonic-boom phenomena, need to approximate solutions for large time. Therefore, they need a good understanding of the behavior of the solutions in these extended regimes in order to be able to simulate them accurately. This issue needs to be treated carefully, as it was already shown in \cite{Pozo:2015}. In that work, the authors proved that a numerical scheme with an acceptable accuracy in short-time intervals could completely disturb the large-time behavior of solutions due to the numerical viscosity introduced by the numerical approximation. It is well known that the asymptotic profile of the inviscid Burgers equation belongs to a two-parameter family of N-waves \cite{Liu:1984}, whereas these N-waves are mere intermediate metastable states of the viscous Burgers equation \cite{Kim:2001}. In our case, \eqref{eq:modelfull} is not a hyperbolic equation and, hence, the asymptotic profile is not an N-wave, but a diffusive wave too. Nevertheless, in our simulations we show that small values for $\nu$ and $c$ require a similar treatment from the numerical point of view, as if the equation was a hyperbolic conservation law. In fact, in those situations, the solution may develop very steep regions (in what follows we refer to these as quasi-shocks), which numerically behave almost like shocks.

Besides the nonlinear term, in this work we also analyze the influence of the operator \eqref{eq:operator} on the large-time behavior of the solutions of the augmented Burgers equation. For the sake of simplicity, we consider only one molecular relaxation phenomenon and homogeneous atmosphere --thus, we neglect the last two terms in \eqref{eq:sonic}. In that framework, note that \eqref{eq:sonic} can be expressed as in \eqref{eq:modelfull}. Moreover, the asymptotical analysis done in the first sections is focused on the case $\nu=c=\theta=1$, but the extension to any positive value of these parameters is immediate. We will omit the subindex $\theta$ whenever its value is one. In this case, we have that
\begin{equation*}
	K*u_{xx}=K*u-u+u_x.
\end{equation*}
Thus, \eqref{eq:modelfull} can be rewritten in a more suitable manner as follows:
\begin{equation}\label{eq:model}
	\begin{cases}
		u_t=uu_x+ u_{xx}+ K*u-u+u_x,&(t,x)\in(0,\infty)\times\R,\\
		u(t=0,x)=u_0(x),&x\in\R.
	\end{cases}
\end{equation}

The main goals of the present paper are to analyze the asymptotic behavior of the solutions to \eqref{eq:model} as $t\to\infty$ and to build a semi-discrete numerical scheme that preserves this behavior. In what concerns the large-time behavior of solutions of system \eqref{eq:model}, the main result is stated in the following theorem.
\begin{thrm}\label{teo:asymptotics}
	Let $u_0\in L^1(\R)$. For any $p\in[1,\infty]$, the solution $u$ to \eqref{eq:model} satisfies
		\begin{equation*}
			t^{\frac12(1-\frac1p)}\|u(t)-u_M(t)\|_p\longrightarrow0,\quad\mbox{as }t\rightarrow\infty,
		\end{equation*}
	where $u_M(t,x)$ is the solution of the following equation:
		\begin{equation*}
			\begin{cases}\label{eq:burgersvisc3}
				u_t=uu_x+2 u_{xx},& x\in\R,t>0,\\
				u(0)=M\delta_0.
			\end{cases}
		\end{equation*}
	Here $\delta_0$ denotes the Dirac measure at the origin and $M$ is the mass of the initial data, $M=\int_\R u_0(x)dx$.
\end{thrm}

In the cases when $\nu$, $c$ and $\theta$ are no longer equal to one, the asymptotic profile does not depend on $\theta$. Moreover, the coefficient in front of the viscosity term in the equation satisfied by the profile is $\nu+c$:
\begin{equation*}
	\begin{cases}
		u_t=uu_x+(\nu + c)u_{xx},& x\in\R,t>0,\\
		u(0)=M\delta_0.
	\end{cases}
\end{equation*}
As a matter of fact, our results are also valid for the case $\nu>0$ and $c=0$, which corresponds to the classical viscous Burgers equation. At the continuous level, this has been long known (see, for instance, \cite{Escobedo:1991} and the references therein). But, to the best of our knowledge, the results for the semi-discrete framework included in our work are new too. On the contrary, the case $\nu=0$ and $c>0$ would require additional results that are beyond the scope of this paper.

Note also that the general case mentioned above will be particularly important at the numerical level. On the one hand, when choosing the numerical flux to discretize the nonlinearity, we need to handle thoroughly the numerical viscosity that is introduced. In \cite{Pozo:2015}, it is shown that in the hyperbolic case, the N-wave asymptotic profile could be destroyed if the numerical flux is not chosen carefully. In our case, if $\nu$ and $c$ are much smaller than $\dx^2/(2\dt)$ ($\dx$ being the mesh-size and $\dt$, the time-step), the Lax-Friedrichs scheme would make the diffusion start dominating much earlier due to the numerical viscosity. On the other hand, we need to treat the truncation of the integral term in such a manner that we do not introduce undesired pathologies in the large-time behavior of the numerical solutions. We do this by means of two correcting factors for the terms $u$ and $u_x$ in \eqref{eq:model}.

Let us denote by $\ud$ an approximation to the solution $u$ of \eqref{eq:model}. We define this piecewise constant in space function as follows:
\begin{equation}\label{eq:udelta}
	\ud(t,x)=u_j(t),\quad x\in(x_{j-1/2},x_{j+1/2}), t\ge0,
\end{equation}
where $x_{j+1/2}=(j+\frac12)\dx$, for all $j\in\Z$, and $\dx>0$ is a given mesh-size. We will also denote by $x_{j}=j\dx$ the intermediate points of the spatial cells. For each $j\in\Z$ we need to compute a function $u_j(t)$ that approximates the value of the solution in the cell. Taking into account the issues enumerated above, we choose the following discretization of \eqref{eq:model}: the Engquist-Osher scheme for the flux, centered finite differences for the laplacian and the composite rectangle rule for the integral:
\begin{equation}\label{eq:scheme}
	\begin{cases}
		u_j'(t)=\dfrac{g_{j+1/2}(t)-g_{j-1/2}(t)}{\dx}+\dfrac{u_{j-1}(t)-2u_j(t)+u_{j+1}(t)}{\dx^2} \\[10pt]
		\quad\qquad\quad + \displaystyle\sum_{m=1}^{N} \omega_m u_{j-m}(t) - \coefzero u_j(t) + \coefone \dfrac{u_{j+1}(t)-u_j(t)}{\dx},&\quad j\in\Z, t\ge0, \\[10pt]
		\displaystyle u_j(0)=\frac{1}{\dx}\int_{x_{j-1/2}}^{x_{j+1/2}} u_0(x)dx, & \quad j\in\Z,
	\end{cases}
\end{equation}
where 
\begin{equation}\label{eq:factorwm}
	\omega_m=\int_{x_{m-1}}^{x_m}K(z)dz=e^{-m\dx} \left(e^{\dx}-1\right),\quad m=1,\dots,N,
\end{equation}
and
\begin{equation*}
	g_{j+1/2}(t)=\frac{u_j(t)\big(u_j(t)-|u_j(t)|\big)}{4}+\frac{u_{j+1}(t)\big(u_{j+1}(t)+|u_{j+1}(t)|\big)}{4}, \quad j\in\Z,t\ge0.
\end{equation*}
The parameter $N=N(\dx)\in\N$ denotes the number of nodes considered in the quadrature formula of the integral. The correcting factors $\coefzero $ and $\coefone $ in front of the approximations of $u$ and $u_x$, given by
\begin{equation}\label{eq:factorfn}
	\coefzero = \sum_{m=1}^N \omega_m \quad \mbox{ and } \quad \coefone = \dx \sum_{m=1}^N m \omega_m,
\end{equation}
handle, from the asymptotic behavior point of view, the correct truncation of the nonlocal term:
\begin{equation*}
	(K\ast u-u+u_x)(x)=\int_{\R} K(x-y)(u(y)-u(x)-(y-x)u_x(x))dy \simeq \sum _{m=1}^N \omega_m\left(u_{j-m}-u_j+m\frac{u_{j+1}-u_j}{\Delta x}\right).
\end{equation*}

Finally, for $\dx$ fixed we study the asymptotic behavior as $t\to\infty$ of these semi-discrete solutions $\ud$.
\begin{thrm}\label{teo:asymptoticsnum}
Let $u_0\in L^1(\R)$, $\dx>0$ and $\ud$ be the corresponding solution of the semi-discrete scheme \eqref{eq:scheme} for the augmented Burgers equation \eqref{eq:model}. For any $p\in[1,\infty]$, the following holds
	\begin{equation}\label{eq:convergencerates}
		t^{\frac12(1-\frac1p)}\|\ud(t)-u_M^\Delta(t)\|_p\longrightarrow0,\quad\mbox{as }t\rightarrow\infty,
	\end{equation}
where $u_M^\Delta(t,x)$ is the unique solution of the following viscous Burgers equation:
	\begin{equation*}
		\begin{cases}
			v_t=vv_x+(1+\coeftwo ) v_{xx},\quad x\in\R,t>0,\\
			v(x,0)=M\delta_0.
		\end{cases}
	\end{equation*}
Here, $M=\int_\R u_0(x)dx$ is the mass of the initial data and
\begin{equation}\label{eq:factorf2}
	\coeftwo =\frac{\dx^2}{2} \left(\sum_{m=1}^N m(m-1)\omega_m \right).
\end{equation}
\end{thrm}

Let us observe that if $N$ is taken such that  $N\dx\to\infty$ when $\dx\to0$ and $N\to\infty$, then $\coeftwo \to1$, which is, precisely, the value that we should expect from the continuous model. Besides, let us remark that in the case where $\nu$, $c$ and $\theta$ are not necessarily equal to one, the asymptotic profile is the unique solution of:
\begin{equation*}
	\begin{cases}
		v_t=vv_x+(\nu+c\,F^{\Delta,\theta}_2) v_{xx},\quad x\in\R,t>0,\\
		v(x,0)=M\delta_0.
	\end{cases}
\end{equation*}
In this case, we take 
\begin{equation*}
	\omega_m^\theta=e^{-m\dx/\theta} \left(e^{\dx/\theta}-1\right)
\end{equation*}
and
\begin{equation*}
	F^{\Delta,\theta}_0= \sum_{m=1}^N \omega_m^\theta, \qquad F^{\Delta,\theta}_1= \frac{\dx}{\theta} \sum_{m=1}^N m \omega_m^\theta \quad\mbox{ and }\quad F^{\Delta,\theta}_2=\frac{\dx^2}{2\theta^2} \left(\sum_{m=1}^N m(m-1)\omega_m^\theta \right).
\end{equation*}
In the same conditions as above, for a fixed $\theta$ we still have that $F^{\Delta,\theta}_2$ converges to one.
 
\begin{rmrk}
We emphasize that function $u_M$ in Theorem \ref{teo:asymptotics} and $u_M^\Delta$ in Theorem \ref{teo:asymptoticsnum} are both particular cases of $u_M^\nu$ ($\nu=2$ and $\nu=1+\coeftwo$ respectively), which is solution of the equation
\begin{equation*}	
	\begin{cases}
		u_t=uu_x+\nu u_{xx},& x\in\R,t>0,\\
		u(0)=M\delta_0.
	\end{cases}
\end{equation*}
In fact, $u_M^\nu$ is explicitly given by (see \cite{Escobedo:1991})
\begin{equation}\label{eq:selfsimburgers}
	u^\nu_M(t,x)=2\sqrt{\nu}\, t^{-\frac12}\exp\left(\frac{-x^2}{4\nu t}\right)\left[C_M+\int_{-\infty}^{x/\sqrt{\nu t}}\exp\left(\frac{-s^2}{4}\right)ds\right]^{-1},
\end{equation}
where $C_M\in\R$ is a constant such that $\int_\R u_M^\nu(t,x)dx=M$, for all $t>0$. This shows that both $u_M$ and $u_M^\Delta$ are of the form $t^{-\frac12}f_M\left(\frac{x}{\sqrt{t}}\right)$ for some function $f_M$ and, hence, self-similar. Note also that $u_M^\Delta\to u_M$ as $\dx\to0$.
\end{rmrk}

Moreover, as we can see in the numerical experiments, the numerical flux needs to be chosen carefully, to avoid adding an extra viscosity term to the equation of the asymptotic profile. This has already been observed in \cite{Pozo:2015} in the context of the numerical approximation of the inviscid Burgers equation. That extra viscosity term, of the order of $\dx^2/(2\dt)$, would affect critically the numerical solution if both parameters $\nu$ and $c$ were much smaller. Note also that taking $\coefzero =\coefone =1$ would add undesired phenomena, such as a transport, to the equation too.

Let us conclude this section by adding a final comment on the time discretization, which we do not address in this paper. At the continuous/semi-discrete level, we obtain estimates on the solution that allow us to prove the compactness of a family of rescaled solutions. Then, the asymptotic behavior is obtained as in \eqref{eq:convergencerates}. The analogous step for the fully time-explicit discrete scheme requires further development.

The paper is organized as follows. In Section \ref{section2}, we deal with the well-posedness of equation \eqref{eq:model} and the asymptotical behavior of its solutions. In Section \ref{section3}, we focus on the semi-discrete numerical scheme \eqref{eq:scheme}, showing its convergence and analyzing for a fixed $\dx$ the large-time behavior of the numerical solutions. To illustrate the main results of this work, we conclude with some numerical simulations in Section \ref{section4}.

In this paper we have considered Engquist-Osher numerical flux for the discretization of the convective term. Let us remark that any other scheme satisfying the analysis done in \cite{Pozo:2015} would be valid too. For instance, Godunov numerical flux would be acceptable, whereas Lax-Friedrichs-type ones are not (as we highlight in Section \ref{section4}).

%%%%%%%%%%%% SECTION 2 %%%%%%%%%%%%%%%%%%%%%%%%%

\section{Analysis of the augmented Burgers equation}\label{section2}

In this section we study the well-posedness of the Cauchy problem for \eqref{eq:model} with initial data in $L^1(\R)$. We also obtain estimates in the $L^p$-norms of its solution, which we subsequently denote $\|\cdot\|_p$. We mainly proceed as in \cite{Escobedo:1991} and \cite{Karch:2010}.

%%%%%%%%%%%% SECTION 2.1 %%%%%%%%%%%%%%%%%%%%%%%%

\subsection{Existence and uniqueness of solutions}\label{section21}
The following theorem concerns the global existence of solutions and specifies their regularity. Let us remark that the result coincides with the one for the classical convection-diffusion equation \cite{Escobedo:1991}.

\begin{thrm}\label{teo:existuni}
For any $u_0\in L^1(\R)$, there exists a unique solution $u\in C([0,\infty), L^1(\R))$ of \eqref{eq:model}. Moreover, it also satisfies 
\begin{equation*}
	u\in C((0,\infty), W^{2,p}(\R))\cap C^1((0,\infty), L^p(\R)),\quad\forall p\in(1,\infty).
\end{equation*}
Additionally, equation \eqref{eq:model} generates a contractive semigroup in $ L^1(\R)$.
\end{thrm}

\begin{proof}
\textit{Existence in $L^1(\R)\cap L^\infty(\R)$.}
The local existence of the solution follows by a classical Banach fixed point argument as in \cite{Escobedo:1991} or \cite{Ignat:2007}. To extend the solution globally, we deduce a priori estimates on the $L^1(\R)$ and $L^\infty(\R)$ norms of the solution. Let us first focus on the $ L^1$-norm. Multiplying \eqref{eq:model} by $\sign(u)$ and integrating in $\R$, it follows that
\begin{equation}\label{eq:estiml1}
	\frac{d}{dt}\int_\R |u|dx\le\int_\R(K*u-u)\sign(u)dx\le \int_\R Kdx \int_\R |u| dx -\int_\R |u| dx \le 0
\end{equation}
and, consequently, $\|u(t)\|_1\le\|u_0\|_1$.

To estimate the $ L^\infty$-norm similar arguments apply. We define $\mu=\|u_0\|_\infty$, multiply equation \eqref{eq:model} by ${\sign[(u-\mu)^+]}$, where $z^+:=\max\{0,z\}$, and integrate it in $\R$. We obtain
\begin{align}\label{eq:estimlinf}
	\frac{d}{dt}\int_\R (u-\mu)^+dx&\le \int_\R (K*u-u+u_x)\sign(u-\mu)^+dx = \int_\R (K*(u-\mu)-(u-\mu))\sign(u-\mu)^+dx\\
	&\leq\int_\R K*(u-\mu)^+ - \int _{\R}(u-\mu)^+\leq 0 \notag .
\end{align}
We conclude that $(u-\mu)^+\le(u_0-\mu)^+=0$ and, consequently, $u(t)\le\mu$ almost everywhere. The same argument for $(u+\mu)^-$, where $z^-:=-\max\{0,-z\}$, shows that $u\ge-\mu$. Therefore, if $u_0\in L^1(\R)\cap L^\infty(\R)$, then $\|u(t)\|_\infty \le \|u_0\|_\infty$ for all $t>0$. Lastly, since both $L^1$-norm and $L^\infty$-norm remain bounded in time, the solution $u$ exists globally.
\\
\textit{Regularity.} It follows from classical regularity arguments (e.g., \cite{Jones:1964}) that
\begin{equation*}
	u\in C((0,T),W^{2,p}(\R))\cap C^1((0,T),L^p(\R))
\end{equation*}
for every $p\in (1,\infty)$. This also holds for $T=\infty$. Let us remark that this regularity makes the integrals in the previous steps well defined.
\\
\textit{Uniqueness.}
To prove the uniqueness of solution it is enough to check that \eqref{eq:model} generates a contractive semigroup in $L^1(\R)$; that is, for any initial datum $u_0,v_0\in L^1(\R)\cap L^\infty(\R)$
\begin{equation}\label{eq:l1contr}
	\|u(t)-v(t)\|_1\le\|u_0-v_0\|_1,\quad\forall t>0,
\end{equation}
where $u$ and $v$ are the corresponding solutions. An analogous argument as in \eqref{eq:estiml1}, applied to the equation verified by $u-v$, shows
\begin{equation*}
	\frac{d}{dt}\int_\R|u-v|dx\le0,
\end{equation*}
hence the contraction property in $L^1(\R)$.
\\
\textit{Existence and uniqueness in $L^1(\R)$.} The extension of the result to a general $u_0\in L^1(\R)$ can be done following the same arguments as in \cite{Escobedo:1991}.
\end{proof}

%%%%%%%%%%%% SECTION 2.2 %%%%%%%%%%%%%%%%%%%%%%%%

\subsection{Decay estimates and large-time behavior}\label{section22}

Now we obtain $L^p$-decay rates for the solution to \eqref{eq:model}. These are the same as the ones for the viscous Burgers equation \cite{Escobedo:1991}.

\begin{prpstn}\label{prop:decay}
For all $p\in[1,\infty]$, there exists a positive constant $C(p)$ such that
	\begin{equation}\label{eq:teodecay}
		\|u(t)\|_p\le C(p) \|u_0\|_1 t^{-\frac12(1-\frac1p)},\quad\forall t>0,
	\end{equation}
for all solutions of equation \eqref{eq:model} with initial data $u_0\in L^1(\R)$.
\end{prpstn}

\begin{proof}
The case $p=1$ is an immediate consequence of Theorem \ref{teo:existuni}. In the case $p\in[2,\infty)$, we multiply equation \eqref{eq:model} by $|u|^{p-2}u$ and integrate it in $\R$. We obtain:
\begin{align}\label{eq:int_simpl}
	\frac1p \frac{d}{dt} \left(\|u\|_p^p\right)=\int_\R|u|^{p-2}uu_tdx&=\int_\R |u|^pu_xdx+\int_\R|u|^{p-2}uu_{xx}dx + \int_\R|u|^{p-2}u (K*u-u+u_x) dx \\
	&= -\frac{4(p-1)}{p^2}\left\|\left(|u|^{p/2}\right)_x\right\|_2^2 -\|u\|_p^p+\int_\R|u|^{p-2}u(K*u) dx \notag.
\end{align}
Let us focus on the last term, so that we can compare it with the $ L^p$-norm of $u$. Young's inequality gives us that
\begin{equation*}
	\left||u(t,x)|^{p-2}u(t,x)u(t,y)\right|=|u(t,x)|^{p-1}|u(t,y)|\le\frac{p-1}{p}|u(t,x)|^p+\frac1p|u(t,y)|^p.
\end{equation*}
Thus, using that $K$ has mass one, it follows:
\begin{equation*}
	\left|\int_\R |u|^{p-2}u(K*u)dx\right| \le \int_\R\int_\R K(x-y)|u(t,x)|^{p-1}|u(t,y)|dydx \le \|u\|_p^p.
\end{equation*}
Plugging this last estimate in \eqref{eq:int_simpl} we have
\begin{equation}\label{eq:est_grad}
	\frac{d}{dt} \left(\|u(t)\|_p^p\right)+\frac{4(p-1)}{p}\left\|\left(|u(t)|^{p/2}\right)_x\right\|_2^2 \le0.
\end{equation}
Finally, with the same arguments as in \cite{Escobedo:1991} we obtain the desired estimate \eqref{eq:teodecay} for any $p\in[2,\infty)$. The case $p=\infty$ follows using the techniques of V\'eron \cite{Veron:1979}. The case $p\in(1,2)$ follows by applying H\"older's inequality and \eqref{eq:teodecay} with $p=1$ and $p=2$.
\end{proof}

Similar estimates can be found for the derivative of the solution of \eqref{eq:model}. Let us define the re-scaled function $u_\lambda$, which will also be used in the following section to obtain the asymptotic profile. For $\lambda>0$ we define
\begin{equation}\label{eq:lambda}
	u_{\lambda}(t,x)=\lambda u(\lambda^2 t,\lambda x).
\end{equation}
The scales are the same as for the Burgers or heat equations. Clearly, $u_\lambda$ is the solution of the following equation:
\begin{equation}\label{eq:modellambda}
	\begin{cases}
		u_{\lambda,t}=u_\lambda u_{\lambda,x}+u_{\lambda,xx}+\lambda^2 (K_{\lambda}*u_\lambda-u_\lambda)+\lambda u_{\lambda,x},&(t,x)\in(0,\infty)\times\R, \\
		u_\lambda(0,x)=u_{\lambda,0}(x)=\lambda u_0(\lambda x),&x\in\R,
	\end{cases}
\end{equation}
where $K_\lambda(z)=\lambda K(\lambda z)$, $z\in\R$.
\begin{prpstn}\label{prop:graddecay}
	For each $p\in[1,\infty]$, there exists a constant $C=C(p,\|u_0\|_1)>0$, such that the solution of equation \eqref{eq:model} satisfies 
		\begin{equation}\label{eq:decayest}
			\|u_x(t)\|_p\le C t^{-\frac12(1-\frac1p)-\frac12},\quad \forall t>0.
		\end{equation}
\end{prpstn}
\begin{proof}
First, note that, for any $\tau>0$,
\begin{equation*}
	\| u_{\lambda,x}(\tau) \|_p = \lambda^{2-\frac{1}{p}} \|u_x(\lambda^2 \tau)\|_p,
\end{equation*}
so proving \eqref{eq:decayest} is equivalent to showing that for some $\tau>0$, $\| u_{\lambda,x}(\tau) \|_p$ is uniformly bounded on $\lambda>0$ and, afterwards, taking $\lambda = \sqrt{t/\tau}$. Let us denote by $D_\lambda^t$ the semigroup associated with the linear problem
\begin{equation*}
	\begin{cases}
		v_t = \lambda^2(K_\lambda*v-v) + \lambda v_x, &(t,x)\in(0,\infty)\times\R,\\
		v(0,x)=v_0(x),&x\in\R.
	\end{cases}
\end{equation*}
It is immediate that $D_\lambda^t$ is non-expansive in $L^p(\R)$, $1\leq p<\infty$,
\begin{equation*}
	\frac 1p\frac{d}{dt}\int_\R |v|^p dx = \lambda^2 \int_\R (K_\lambda*v-v) |v|^{p-1}\sign(v) dx \le 0.
\end{equation*}
On the other hand, for all $\tau>0$, function $u_\lambda$ solution of \eqref{eq:modellambda} verifies the following integral equation:
\begin{align*}
	u_\lambda(t+\tau)&=G(t)*D_\lambda^{t} u_\lambda(\tau)+\int_0^t G(t-s)*D_\lambda^{t-s} \left(\left(\frac{u_\lambda^2(s+\tau)}{2}\right)_x\right)ds,
\end{align*}
where $G(t)$ is the heat kernel. Using the fact that $D_\lambda$ is non-expansive in $L^p(\R)$, $1\leq p<\infty$, and following the same arguments as in \cite{Escobedo:1991} we obtain the desired results. For complete details see \cite{Pozo:2014}.
\end{proof}

%%%%%%%%%%%% SECTION 2.3 %%%%%%%%%%%%%%%%%%%%%%%%

\subsection{Asymptotic expansion}\label{section23}

The decay rates of the previous section will allow us to obtain the asymptotic profile of solutions for \eqref{eq:model}. Note that taking $\lambda=\sqrt{t}$ we have that
\begin{equation*}
	\|u_\lambda(1)-u_M(1)\|_1=\|u(t)-u_M(t)\|_1,
\end{equation*}
due to the definition of $u_\lambda$ and the self-similar nature of $u_M$. Thus, the aim is to compute the limit $\lambda\rightarrow\infty$ in \eqref{eq:modellambda}, which is equivalent to taking the limit $t\to\infty$ in \eqref{eq:model} when $p=1$.

Let us first observe that the estimates in Proposition \ref{prop:decay} and Proposition \ref{prop:graddecay} are also valid for $u_\lambda$ defined in \eqref{eq:lambda}. The mass is conserved too. We state this in the following lemma.
\begin{lmm}\label{lemma:decayulambda}
	For each $p\in[1,\infty]$, there exists a constant $C=C(p,\|u_0\|_1)>0$ such that, for all $\lambda>0$, the solution of \eqref{eq:modellambda} satisfies 
		\begin{equation*}
			\|u_\lambda(t)\|_p\le C t^{-\frac12(1-\frac1p)} \qquad \mbox{ and } \qquad \|u_{\lambda,x}(t)\|_p\le C t^{-\frac12(1-\frac1p)-\frac12},\,\quad\forall t>0.
		\end{equation*}
	Moreover, the mass of $u_\lambda$ is conserved in time.
\end{lmm}

\begin{proof}
We just have to use the definition of $u_\lambda$ in \eqref{eq:lambda} and apply Proposition \ref{prop:decay}. For all $t>0$ and $\lambda>0$ we have
\begin{equation*}
	\|u_\lambda(t)\|_p=\lambda^{1-\frac1p}\|u(\lambda^2t)\|_p \le Ct^{-\frac12(1-\frac1p)}.
\end{equation*}
Same procedure applies to $u_{\lambda,x}$, concerning Proposition \ref{prop:graddecay}. Regarding the last result, it is easy to see that:
\begin{equation*}
	\int_\R u_\lambda(t,x)dx=\int_\R u(\lambda^2t,x)dx=\int_\R u_0(x)dx,
\end{equation*}
which proves the mass conservation.
\end{proof}

In particular, this lemma implies that, for any finite time interval $[\tau,T]$ with $0<\tau<T<\infty$, the set $\{u_\lambda\}_{\lambda>0}$ is uniformly bounded in $ L^\infty([\tau,T], L^p(\R))$, $1\le p \le \infty$.

%%%%%%%%%%%% SECTION 2.3.1 %%%%%%%%%%%%%%%%%%%%%%%%

\subsubsection{Compactness of the family $\{u_\lambda\}_{\lambda>0}$}

As we said at the beginning, we would like to pass to the limit $\lambda\rightarrow\infty$. We need the following theorem due to J. Simon \cite{Simon:1987}, as an extension of the Aubin-Lions Lemma, to assure the compactness of the set $\{u_\lambda\}_{\lambda>0}$.

\begin{thrm}[{\cite[Theorem 5]{Simon:1987}}]\label{teo:simon}
	Let $X$, $Z$ and $Y$ be Banach spaces satisfying $X\subset Z\subset Y$ with compact embedding $X\subset Z$. Assume, for $p\in[1,\infty]$ and $T>0$, that $F$ is bounded in $ L^p(0,T;X)$ and $\{\partial_tf:f\in F\}$ is bounded in $ L^p(0,T;Y)$. Then, $F$ is relatively compact in $ L^p(0,T;Z)$ and, in the case of $p=\infty$, also in $C(0,T;Z)$.
\end{thrm}

Applying this result we can prove the following theorem regarding the relative compactness of the set $\{u_\lambda\}_{\lambda>0}$. In the sequel, for any functions $f$ and $g$, we denote $f \lesssim g$ if there exists a constant $C>0$, not depending on the scaling parameter nor the time, such that $f\le Cg$.

\begin{thrm}\label{teo:compact}
	For every $0<\tau<T<\infty$, the set $\{u_\lambda\}_{\lambda>0}\subset C([\tau,T], L^1(\R))$ is relatively compact.
\end{thrm}
\begin{proof}
\textit{Step 1: Compactness in  $C([\tau,T], L^1_{loc}(\R))$.} First, for any $r>0$ we will show the relative compactness in $C([\tau,T], L^2(I))$, with $I=[-r,r]$. Let us consider the spaces $X= H^1(I)$, $Z= L^2(I)$ and $Y= H^{-1}(I)$. We would like to apply Theorem \ref{teo:simon} to the set $F=\{u_\lambda\}_{\lambda>0}$.

From Lemma \ref{lemma:decayulambda} we know that $\{u_\lambda\}_{\lambda>0}$ and $\{u_{\lambda,x}\}_{\lambda>0}$ are bounded in $L^\infty([\tau,T],L^2(I))$. In particular, the first condition of Theorem \ref{teo:simon} on $F$ is fulfilled. Therefore, it suffices to check that $u_{\lambda,t}$ is bounded in $ L^\infty([\tau,T], H^{-1}(I))$. Using \eqref{eq:modellambda}, for every $\varphi\in C^\infty_c(I)$, we have:
	\begin{align} \label{eq:localcomp}
		\left|\int_\R u_{\lambda,t}\varphi dx\right|&\le\left|\int_\R u_\lambda\ u_{\lambda,x}\varphi dx\right|+ \left|\int_\R u_{\lambda,xx}\varphi dx\right|+\left|\int_\R \big(\lambda^2(K_{\lambda}*u_{\lambda}-u_\lambda)+\lambda u_{\lambda,x}\big)\varphi dx\right| \\
		&\lesssim\|\varphi_x\|_2\|u_\lambda\|^2_4+\|\varphi_x\|_2\|u_{\lambda,x}\|_2+\left|\int_\R \big(\lambda^2(K_{\lambda}*u_{\lambda}-u_\lambda)+\lambda u_{\lambda,x}\big)\varphi dx\right|. \notag
	\end{align}
Obviously, the first two terms on the right hand side of \eqref{eq:localcomp} are uniformly bounded in $[\tau,T]$, so let us focus on the third one:
\begin{equation*}
	\mathcal{I}_{\lambda}=\left|\int_\R \big(\lambda^2(K_{\lambda}*u_{\lambda}-u_\lambda)+\lambda u_{\lambda,x}\big)\varphi dx\right| =\left|\int_\R \Big(\lambda^2\big(\widehat{K}(\xi/\lambda) -1\big)+i\lambda\xi\Big)\widehat{u_{\lambda}}(\xi)\widehat{\varphi}(\xi) d\xi\right|.
\end{equation*}
Let us denote
\begin{equation*}
	m_{\lambda}(\xi)=\lambda^2\Big(\widehat{K}(\frac{\xi}{\lambda}) -1\Big)+i \lambda\xi.
\end{equation*}
We claim that
\begin{equation}\label{eq:mfourier}
	|m_\lambda(\xi)|\le \xi^2,\quad \forall \xi\in\R,\forall \lambda>0.
\end{equation}
Using the Cauchy-Schwartz inequality, we have:
\begin{equation}\label{eq:estimcomp}
	\mathcal{I}_{\lambda}=\left|\int_\R m_{\lambda}(\xi)\widehat{u_{\lambda}}(\xi)\widehat{\varphi}(\xi) d\xi\right|\lesssim \left\|\varphi\right\|_{ H^1(\R)}\left\|u_\lambda\right\|_{ H^1(\R)}.
\end{equation}
Hence, going back to \eqref{eq:localcomp} and replacing \eqref{eq:estimcomp}, we obtain
\begin{equation*}
	\left|\int_\R u_{\lambda,t}\varphi dx\right|\lesssim \left\|\varphi\right\|_{ H^1(\R)}\left(\|u_\lambda\|^2_4+\left\|u_\lambda\right\|_{ H^1(\R)}\right).
\end{equation*}
By Lemma \ref{lemma:decayulambda}, all the quantities in the right-hand side are uniformly bounded in $[\tau, T]$. Consequently, the set $\{u_\lambda\}_{\lambda>0}$ is relatively compact in $C([\tau,T], L^2(I))$.

It remains to prove claim \eqref{eq:mfourier}. Observe that
\begin{equation}\label{eq:fourierK}
	\widehat{K}(\xi)=\frac{1}{1+i\xi},\quad\xi\in\R,
\end{equation}
and, therefore,
\begin{equation*}
	|m_\lambda(\xi)|=\left|\lambda^2\left(\frac{1}{1+i\xi/\lambda} -1\right)+i\lambda\xi\right|=\frac{\lambda\xi^2}{\sqrt{\lambda^2+\xi^2}}\le\xi^2,\quad \forall\lambda>0.
\end{equation*}
Since $L^2(I)$ is continuously embedded in $L^1(I)$, the compactness in $C([\tau,T], L^2(I))$ is clearly transferred to $C([\tau,T], L^1(I))$. To extend this local result to the globally we prove uniform, with respect to $\lambda$, estimates on the tails of $u_\lambda$.

\textit{Step 2: Uniform control of the tails.} For every $r>0$, let us define function $\psi_r(z)=\psi(z/r)$, where $\psi$ is a nonnegative $C^\infty(\R)$ function such that 
\begin{equation}\label{eq:psi}
	\psi(z)=
	\begin{cases}
		0,&|z|<1,\\ 1,&|z|>2.
	\end{cases}
\end{equation}
Since $\{u_\lambda\}_{\lambda>0}$ is relatively compact in $C([\tau,T],L^1(I))$, it suffices to show that
\begin{equation}\label{eq:sup}
	\sup_{t\in [\tau,T]}\|u_\lambda(t)\psi_r\|_1\longrightarrow0 \quad\mbox{as }r\rightarrow\infty,\mbox{ uniformly for } \lambda>0.
\end{equation}

We first observe that it is enough to consider nonnegative initial data. For any $u_0,v_0\in L^1(\R)$ satisfying $u_0\leq v_0$ we can show that the corresponding solutions of \eqref{eq:modellambda} satisfy $u_\lambda\leq v_\lambda$. Thus choosing $v_0=|u_0|$ and $v_0=-|u_0|$ we can show that $|u_\lambda|\leq \tilde u_\lambda$ where $\tilde u_\lambda$ is the solution corresponding to $|u_0|$ initial data. This reduces \eqref{eq:sup} to the case of non-negative solutions. Let us assume that $u_\lambda$ is a nonnegative solution. We multiply \eqref{eq:modellambda} by $\psi_r$ and integrate it over $(0,t)\times\R$. We obtain:
\begin{align*}
	\int_0^t\int_\R u_{\lambda,s}\psi_r dxds&=-\frac12 \int_0^t\int_\R u_\lambda^2 \psi'_r dxds + \int_0^t\int_\R u_\lambda \psi''_r dxds\\
	&\quad\quad+\int_0^t\int_\R \Big(\lambda^2(K_{\lambda}*u_{\lambda}- u_{\lambda})+\lambda u_{\lambda,x}\Big)\psi_{r}dxds.
\end{align*}
and, therefore,
\begin{align}\label{eq:int_j}
	\int_\R u_{\lambda}(t)\psi_r dx&\le \int_\R u_{\lambda,0} \psi_r dx+\frac{\|\psi'\|_\infty}{2r}\int_0^t \|u_\lambda(s)\|_2^2ds+\frac{\|\psi''\|_\infty}{r^2}\int_0^t\|u_{\lambda}(s)\|_1ds \\
	&\quad+\int_0^t \int_\R \Big(\lambda^2\big(K_{\lambda}*u_{\lambda}(s)- u_{\lambda}(s)\big)+ \lambda u_{\lambda,x}(s)\Big)\psi_rdx ds. \notag
\end{align}
We have to obtain an estimate on the last term in the integral, uniformly on $\lambda$. Let us denote
\begin{equation*}
	\mathcal{J}=\int_\R \Big(\lambda^2\big(K_{\lambda}*u_{\lambda}(s)- u_{\lambda}(s)\big)+ \lambda u_{\lambda,x}(s)\Big)\psi_rdx.
\end{equation*}
A change of variables and integration by parts give us that
\begin{align}\label{eq:estimJ}
	\mathcal{J}&=\lambda^2 \int _{\R}\int _{\R}K(x-y)u(\lambda^2s,y)\psi_{\lambda r}(x)dydx - \lambda^2 \int _{\R}u(\lambda^2s,x)\psi_{\lambda r}(x)dx-\lambda^2\int _{\R} u(\lambda^2s,x) \psi_{\lambda r}'(x)dx \\
	& = \lambda^2\int _{\R} u(\lambda^2s,y) \Big(\int_{\R} K(x-y)\big(\psi_{\lambda r}(x)-\psi _{\lambda r}(y)-(x-y)\psi_{\lambda r}'(y)\big)dx\Big)dy \notag \\
	&\leq \lambda^2 \|u_0\|_{1} \|(\psi _{\lambda r})''\|_{\infty} = \frac{\|u_0\|_{1} \|\psi''\|_{\infty}}{r^2}. \notag
\end{align}
Plugging \eqref{eq:estimJ} into \eqref{eq:int_j} and using Proposition \ref{prop:decay}, we get:
\begin{equation*}
	\int_\R u_\lambda(t) \psi_r dx \le \int_\R u_{0} \psi_{\lambda r} dx+C\left(\frac{\sqrt{t}}{r}+\frac{t}{r^2}\right)
\end{equation*}
where $C>0$ depends only on $\|u_0\|_1$ and $\|\psi\|_{W^{2,\infty}(\R)}$, which are both bounded. For $\lambda>1$, since $\psi_{r}(x)>\psi_{\lambda r}(x)$, we get
\begin{equation*}
	\int_\R u_\lambda(t,x) \psi_r(x) dx \le \int_\R u_{0}(x) \psi_{r}(x) dx+C\left(\frac{\sqrt{t}}{r}+\frac{t}{r^2}\right),
\end{equation*}	
which tends to zero uniformly on $\lambda$ when $r\rightarrow\infty$. Therefore, we proved \eqref{eq:sup} and, consequently, we can assure that $\{u_\lambda\}_{\lambda>0}$ is relatively compact in $ C([\tau,T], L^1(\R))$.
\end{proof}

Modifying slightly the previous proof, we can also conclude the following lemma, regarding the initial condition $u_{\lambda,0}$.
\begin{lmm}\label{lemma:decayulambda3}
	For every test function $\varphi\in C^2_b(\R)$, there exists a constant $C=C(\varphi,u_0)>0$, such that
	\begin{equation*}
		\left|\int_\R u_\lambda(t,x)\varphi(x)dx-\int_\R u_{\lambda,0}(x)\varphi(x)dx\right|\le C(t+\sqrt{t}),\quad\forall t>0,
	\end{equation*}
	holds uniformly on $\lambda>0$.
\end{lmm}

\begin{proof}
We multiply \eqref{eq:modellambda} by $\varphi\in C^2_b(\R)$ and integrate it over $(0,t)\times\R$. We get:
\begin{align*}
	\int_0^t\int_\R u_{\lambda,t}\varphi &=\int_0^t\int_\R u_\lambda u_{\lambda,x}\varphi +\int_0^t\int_\R u_{\lambda,xx}\varphi +\int_0^t\int_\R \big(\lambda^2(K_{\lambda}*u_{\lambda}-u_{\lambda})+\lambda u_{\lambda,x}\big)\varphi .
\end{align*}
Integrating by parts and making use of Lemma \ref{lemma:decayulambda}, we have
\begin{align*}		
	\left|\int_\R u_{\lambda}(t)\varphi dx-\int_\R u_{\lambda,0}\varphi dx\right| &\le \frac{\|\varphi'\|_\infty}{2} \int_0^t \|u_\lambda(s)\|_2^2ds+\|\varphi''\|_\infty\int_0^t\|u_{\lambda}(s)\|_1ds\\
	&\quad\quad+\left|\int_0^t\int_\R \Big(\lambda^2\big(K_{\lambda}*u_{\lambda}(s)-u_{\lambda}(s)\big)+\lambda u_{\lambda,x}(s)\Big)\varphi dxds\right|.
\end{align*}
To conclude the proof, it is enough to apply a similar argument as for \eqref{eq:int_j} to get:
\begin{align*}		
	\left|\int_\R u_\lambda(t)\varphi dx-\int_\R u_{\lambda,0}\varphi dx\right| \le C(\|\varphi\|_{W^{2,\infty}(\R)},\|u_0\|_1)(\sqrt{t} +t).
\end{align*}
The proof is now finished.
\end{proof}

%%%%%%%%%%%% SECTION 2.3.2 %%%%%%%%%%%%%%%%%%%%%%%%

\subsubsection{Passing to the limit}

Now we have all the ingredients that we need to prove our main result on the large-time behavior of solutions to problem \eqref{eq:model}, stated in Theorem \ref{teo:asymptotics}.

\begin{proof}[Proof of Theorem \ref{teo:asymptotics}]
	By Theorem \ref{teo:compact}, we know that for every $0<\tau<T<\infty$, the family $\{u_\lambda\}_{\lambda>0}$ is relatively compact in $ C([\tau,T], L^1(\R))$. Consequently, there exists a subsequence of it (which we will not relabel) and a function $\bar u\in C((0,\infty), L^1(\R))$ such that
\begin{equation}\label{eq:sequence}
	u_{\lambda}\longrightarrow \bar u\in C([\tau,T], L^1(\R)),\quad\mbox{as }\lambda\rightarrow\infty.
\end{equation}
We can also assume that $u_{\lambda}(t,x)\rightarrow\bar u(t,x)$ almost everywhere in $(0,\infty)\times\R$ as ${\lambda\rightarrow\infty}$.

Our claim is that, passing to the limit $\lambda\rightarrow\infty$, we obtain that $\bar u$ is a weak solution of the equation:
\begin{equation}\label{eq:burgersvisc}
	\begin{cases}
		\bar u_{t}=\bar u \bar u_{x}+2\bar u_{xx},&(t,x)\in(0,\infty)\times\R,\\
		\bar u(0)=M\delta_0.
	\end{cases}
\end{equation}
Let us multiply equation \eqref{eq:modellambda} by a test function $\phi\in C^\infty_c((0,\infty)\times\R)$ and integrate it over $(0,\infty)\times\R$. We have:
\begin{equation*}
	-\int_0^\infty \int_\R u_{\lambda} \phi_t = \int_0^\infty \int_\R u_\lambda\ u_{\lambda,x} \phi + \int_0^\infty \int_\R u_{\lambda,xx} \phi+ \int_0^\infty \int_\R (\lambda^2 (K_{\lambda}*u_\lambda-u_\lambda)+\lambda u_{\lambda,x}) \phi.
\end{equation*}
Using the properties of $\{u_\lambda\}_{\lambda>0}$ shown in the previous section, it is sufficient to check that
\begin{equation*}
	\lim_{\lambda\rightarrow\infty} \int_0^\infty \int_\R\Big(\lambda^2 \big(K_{\lambda}*u_\lambda(t)- u_\lambda(t)\big)+ \lambda u_{\lambda,x}(t)\Big)\phi (t)dxdt= \int_0^\infty \int_\R\bar u(t)\phi_{xx}(t)dxdt.
\end{equation*}
Let us focus on the integral over the spatial domain. Taking into account the definition of $K_\lambda$ and that $\int_\R z^m K(z)dz=m! $ for $m\in\N\cup\{0\}$ we obtain
\begin{align}\label{eq:convlap}
	\mathcal{L_\lambda}(t)&:=\int_\R\Big(\lambda^2 \big(K_{\lambda}*u_\lambda(t)- u_\lambda(t)\big)+ \lambda u_{\lambda,x}(t)\Big)\phi (t) dx\\
	&=\int_\R\int_\R u_\lambda(t,x)K(y) \left(\phi(t,x+\frac{y}{\lambda})-\phi(t,x)-\frac{y}{\lambda}\phi_x(t,x)\right)dxdy  \notag\\
	&=\lambda^2 \int_\R \int_\R \Big(\phi(t,x+\frac{y}{\lambda})-\phi(t,x)\Big)K(y)u_\lambda(t,x)dydx -\lambda\int_\R\int_\R \phi_x(t,x)yK(y) u_{\lambda}(t,x)dydx, \notag
\end{align}
Now, because of Taylor's Theorem, we know that there exists a point $\zeta\in(x,x+y/\lambda)$ such that
\begin{equation*}
	\phi(t,x+\frac{y}{\lambda})-\phi(t,x)=\frac{y}{\lambda}\phi_x(t,x)+\frac{1}{2}\frac{y^2}{\lambda^2}\phi_{xx}(t,x)+\frac16\frac{y^3}{\lambda^3}\phi_{xxx}(t,\zeta).
\end{equation*}
We introduce this in \eqref{eq:convlap}:
\begin{align}\label{eq:convfin}
	 \mathcal{L_\lambda}(t) &=\frac{1}{2}\int_\R \phi_{xx}(t,x)u_\lambda(t,x)dx\int_\R y^2K (y) dy + O(\|\phi_{xxx}(t)\|_\infty) \frac{1}{6\lambda}\int_\R u_\lambda(t,x)dx \int_\R y^3K (y) dy\notag \\
	 &=\int_\R \phi_{xx}(t,x)u_\lambda(t,x)dx+O(\|\phi_{xxx}(t)\|_\infty) \frac{1}{\lambda} \int_\R u_\lambda(t,x)dx \notag \\
	 &=\int_\R \phi_{xx}(t,x)u_\lambda(t,x)dx+ \lambda^{-1} O(\|\phi_{xxx}(t)\|_\infty). \notag 
\end{align}
Since $u_\lambda(t)\rightarrow\bar u(t)$ in $C([\tau,T],L^1(\R))$ and $\phi$ has compact support, we obtain that
\begin{equation*}
	\lim_{\lambda\rightarrow\infty} \int_0^\infty \mathcal{L_\lambda}(t)dt= \int_0^\infty \int_\R \phi_{xx}(t,x)\bar u(t,x)dxdt.
\end{equation*}
It follows that $\bar u$ satisfies
\begin{equation*}
	-\int_0^\infty \int_\R \bar u \phi_t = -\frac12 \int_0^\infty \int_\R \bar u^2 \phi_x + 2 \int_0^\infty \int_\R \bar u \phi _{xx}.
\end{equation*}
It remains to identify the behavior of $\bar u$ as $t\to0$. From Lemma \ref{lemma:decayulambda3}, for any $\varphi\in C^2_b(\R)$ we have
\begin{equation*}
	\left|\int_\R u_\lambda(t,x)\varphi(x)dx - \int_\R u_{\lambda,0}(x)\varphi(x)dx \right| \le C(t+\sqrt{t})
\end{equation*}
and, due to \eqref{eq:sequence} and the definition of $u_\lambda$ in \eqref{eq:modellambda}, we deduce by letting $\lambda\to\infty$ that 
\begin{equation*}
	 \lim _{t\downarrow 0} \int_\R \bar u(t,x)\varphi(x)dx = M\varphi(0).
\end{equation*}
Using classical approximation arguments together with the uniform tail control of $u_\lambda$ in \eqref{eq:sup}, we conclude that $\bar u(0)=M\delta_0$ in the sense of bounded measures.

Therefore, we can finally conclude that $\bar u$ is the unique solution $u_M$ of \eqref{eq:burgersvisc}, and that, indeed, the whole family $\{u_\lambda\}_{\lambda>0}$ converges to $u_M$ in $ C((0,\infty), L^1(\R))$. In particular, we have:
\begin{equation*}
	\lim_{\lambda\rightarrow\infty}\|u_{\lambda}(1)-u_M(1)\|_1=0.
\end{equation*}
Setting $\lambda=\sqrt{t}$ and using the self-similar form of $u_M$ (see e.g. \cite{Escobedo:1991}), we obtain that
\begin{equation}\label{eq:limnorm1}
	\lim_{t\rightarrow\infty}\|u(t)-u_M(t)\|_1=0.
\end{equation}

Finally, the convergence in the $ L^p$-norms for $p\in(1,\infty)$ follows from \eqref{eq:limnorm1}, the decay estimate given in Lemma \ref{lemma:decayulambda} for $p=\infty$ and the H\"older inequality. In fact, we have:
\begin{equation}\label{eq:limnormp}
	\|u(t)-u_M(t)\|_p\le\left(\|u(t)\|_\infty+\|u_M(t)\|_\infty\right)^{1-\frac1p}\|u(t)-u_M(t)\|_1^{\frac1p}\le o(t^{-\frac12(1-\frac1p)}).
\end{equation}
In the case of the $ L^\infty$-norm, we use the decay of $u_x(t)$ given by Proposition \ref{prop:graddecay} and the estimate $\|u_{M,x}(t)\|_2\lesssim t^{-\frac34}$, resulting from the explicit formula \eqref{eq:selfsimburgers}. Using the Gagliardo-Nirenberg interpolation inequality and \eqref{eq:limnormp}, we obtain:
\begin{equation}
	\|u(t)-u_M(t)\|_\infty\lesssim \left(\|u_x(t)\|_2+\|u_{M,x}(t)\|_2\right)^{\frac12}\|u(t)-u_M(t)\|_2^{\frac12}\le o(t^{-\frac12}).
\end{equation}
The proof is now finished.
\end{proof}

%%%%%%%%%%%% SECTION 3 %%%%%%%%%%%%%%%%%%%%%%%%%

\section{Semidiscrete scheme}\label{section3}
In this section, we focus on the semi-discrete numerical scheme for equation \eqref{eq:model}, defined in \eqref{eq:scheme}. In order to prove Theorem \ref{teo:asymptoticsnum}, we need some preliminary results on the decay of $u_\Delta$ similar to those obtained in Section \ref{section2} for the solution of equation \eqref{eq:model}. For simplicity, for every $h>0$, we define the operators $d^+_h$ and $d^-_h$ as follows: 
\begin{equation*}
	d^+_hf(x):=\frac{f(x+h)-f(x)}{h} \quad\mbox{ and }\quad d^-_hf(x):=\frac{f(x)-f(x-h)}{h}.
\end{equation*}

As in the continuous case, for $\mu>0$ we also introduce the family of rescaled solutions
\begin{equation}\label{eq:umudef}
	\umu(t,x)=\mu \ud(\mu^2t,\mu x),\quad t\ge0,x\in\R,
\end{equation}
and analyze the behavior of $\umu$ when $\mu\to\infty$. Note that function $\umu$ is piecewise constant on space intervals of length $\dx/\mu$. Moreover, it satisfies the following system:
\begin{equation}\label{eq:umu}
	\begin{cases}
		\displaystyle \umu_t(t,x)= \frac14 \Big(d^+_{\dx/\mu} \big(\umu(t,x)^2\big)+d^-_{\dx/\mu} \big(\umu(t,x)^2\big)\Big) \\[10pt]
		\qquad\qquad\quad + \dx\, d^+_{\dx/\mu} R\big(\umu(t,x-\frac{\dx}{\mu}),\umu(t,x)\big) + d^-_{\dx/\mu} \left(d^+_{\dx/\mu} \umu(t,x) \right) \\[5pt]
		\qquad\qquad\quad +\mu^2{\displaystyle \sum_{m=1}^{N}} \omega_m \umu(t,x-m\frac{\dx}{\mu}) -\mu^2 \coefzero \umu(t,x) +\mu \coefone d^+_{\dx/\mu} \umu(t,x),& t>0,\mbox{ a.e. }x\in\R, \\[10pt]
		\umu(0,x)=\mu\ud^0(\mu x), & \mbox{a.e. }x\in\R,
	\end{cases}
\end{equation}
where
\begin{equation}\label{eq:numviscR}
	R(u,v)=\frac{1}{4\dx}(v|v|-u|u|).
\end{equation}
Of course, the approximated solution $\ud$ defined in \eqref{eq:udelta} and \eqref{eq:scheme} satisfies \eqref{eq:umu} when $\mu=1$. Let us recall that this decomposition of the numerical flux is called the viscous form of the scheme (see, for instance, Chapter III in \cite{Godlewski:1991}). Note also that $R$ is homogeneous of degree 2, allowing the term containing it in \eqref{eq:umu} to disappear as $\mu\to\infty$, as observed in \cite{Pozo:2015}.

For any initial data $\ud^0\in L^1(\R)$, there exists a unique solution in $C^1([0,\infty),L^1(\R))$ for \eqref{eq:umu}. The local existence is obtained by  Banach's fixed point argument, whereas the following Lemma \ref{lemma:conservation}  excludes blow-ups.
Let us remark that the solution $\umu$ of system \eqref{eq:umu} conserves the mass of the initial data $\ud^0$. In fact, note that it is the same as the mass of $u_0$, when $\ud^0$ is defined as in \eqref{eq:scheme}. Moreover, we show that \eqref{eq:umu} defines a contractive semigroup. This will be useful to obtain the estimates for the compactness of $\{\umu\}_{\mu>0}$. For the sake of clarity, we prove this lemma in the Appendix \ref{appendix}.

\begin{lmm}\label{lemma:conservation}
	For any initial data $\ud^0\in L^1(\R)$, the solution $\umu$ to \eqref{eq:umu} satisfies
	\begin{equation*}
		\int_\R \umu (t,x)=\int_\R \ud^0(x),\quad\forall t>0.
	\end{equation*}
	Moreover, \eqref{eq:umu} defines a contractive semigroup in $L^1(\R)$.
\end{lmm}

%%%%%%%%%%% SUBSECTION 3.1 %%%%%%%%

\subsection{$L^1$-$L^p$ estimates}
We are interested in the large-time behavior of $\ud$. The following two propositions are the discrete versions of Proposition \ref{prop:decay} and Proposition \ref{prop:graddecay}. The way of proceeding will be, indeed, very similar.

\begin{prpstn}\label{prop:est1}
For all $p\in[1,\infty]$ there exists a positive constant $C(p)$ such that:
\begin{equation}
	\|\umu(t)\|_p \le C(p)\|\ud^0\|_1 t^{-\frac12(1-\frac1p)},\quad \forall t>0.
\end{equation}
for all solutions of \eqref{eq:umu} with initial data $\ud^0\in L^1(\R)$.
\end{prpstn}
\begin{proof}
The case $p=1$ follows from Lemma \ref{lemma:conservation} with $C(p)=1$. Let us consider the case $\mu=1$ and $p\in[2,\infty)$. We multiply \eqref{eq:umu} by $|\ud|^{p-2}\ud$ and integrate it over the whole space domain. We have:
\begin{align}\label{eq:ineq1112}
	\frac1p \frac{d}{dt} \|\ud(t)\|_p^p&\le I_1 + \int_\R d^-_{\dx}\left(d^+_{\dx}\ud(t,x)\right) |\ud(t,x)|^{p-2}\ud(t,x)dx + I_2,
\end{align}
where
\begin{align*}
	I_1&=\frac14 \int_\R \Big(d^+_{\dx} \big(\ud(t,x)^2\big)+d^-_{\dx} \big(\ud(t,x)^2\big)\Big) |\ud(t,x)|^{p-2}\ud(t,x)dx \\
	&\qquad+ \dx \int_\R d^+_{\dx} R\big(\ud(t,x-\dx),\ud(t,x)\big) |\ud(t,x)|^{p-2}\ud(t,x)dx
\end{align*}
and
\begin{align*}
	I_2&=\int_\R \left(\sum_{m=1}^{N} \omega_m \ud(t,x-m\dx) - \coefzero \ud(t,x) + \coefone d^+_{\dx} \ud(t,x)\right)|\ud(t,x)|^{p-2}\ud(t,x) dx \\
	&=\sum_{m=1}^{N} \omega_m \left(\int_\R \ud(t,x-m\dx) |\ud(t,x)|^{p-2}\ud(t,x) dx - \int_\R |\ud(t,x)|^{p} dx \right) \\
	& \quad + \frac{\coefone }{\dx} \left(\int_\R \ud(t,x+\dx) |\ud(t,x)|^{p-2}\ud(t,x) dx - \int_\R |\ud(t,x)|^{p} dx \right).
\end{align*}
On the following, we will not make explicit the time dependence unless this is necessary.

Now, on the one hand, for any $k\in\Z$, we have that
\begin{equation*}
	\int_\R \ud(x+k \dx)|\ud(x)|^{p-2} \ud(x) dx \le \frac{p-1}{p}\int_\R |\ud(x+k\dx)|^p dx +\frac1p \int_\R |\ud(x)|^p dx = \int_\R |\ud(x)|^p dx.
\end{equation*}
Therefore, $I_2\le0$. 

On the other hand, for $i\in\{-1,0,1\}$ let us denote
\begin{equation*}
	U^\pm_i=\{x\in\R:\pm \ud(x+i\dx)>0\} \quad \text{and} \quad U^0_i=\{x\in\R:\ud(x+i\dx)=0\}.
\end{equation*}	
From the definition of $R$ in \eqref{eq:numviscR}, reordering $I_1$ we get:
\begin{align*}
	I_1& = \frac{1}{4\dx}\int_\R \left(\ud^2(x+\dx)+\ud(x+\dx)|\ud(x+\dx)|\right) |\ud(x)|^{p-2}\ud(x)dx - \frac{1}{2\dx}\int_\R |\ud(x)|^{p+1}dx \\ 
	& \quad + \frac{1}{4\dx}\int_\R \left(\ud(x-\dx)|\ud(x-\dx)|-\ud^2(x-\dx)\right) |\ud(x)|^{p-2}\ud(x)dx \\
	& \le \frac{1}{2\dx}\int_{U_0^+\cap U^+_1} \ud^2(x+\dx) |\ud(x)|^{p-1}dx - \frac{1}{2\dx}\int_\R |\ud(x)|^{p+1}dx \\ 
	& \quad + \frac{1}{2\dx}\int_{U_{-1}^-\cap U_0^-} \ud^2(x-\dx) |\ud(x)|^{p-1}dx.
\end{align*}
Using the following inequality
\begin{equation*}
	a^2|b|^{p-1} \le \frac{2}{p+1}|a|^{p+1} + \frac{p-1}{p+1}|b|^{p+1},\quad\forall a,b\in\R,
\end{equation*}
we obtain that 
\begin{align*}
	I_1 & \le \frac{1}{2\dx} \left(\frac{2}{p+1} \int_{U_{0}^+\cap U_1^+} |\ud(x+\dx)|^{p+1}dx + \frac{p-1}{p+1} \int_{U_0^+\cap U_1^+} |\ud(x)|^{p+1}dx \right) - \frac{1}{2\dx}\int_\R |\ud(x)|^{p+1}dx \\ 
	& \quad + \frac{1}{2\dx} \left(\frac{2}{p+1} \int_{U_{-1}^-\cap U_0^-} |\ud(x-\dx)|^{p+1}dx + \frac{p-1}{p+1} \int_{U_{-1}^-\cap U_0^-} |\ud(x)|^{p+1}dx\right) \\
	&=\frac{1}{2\dx} \left(\frac{2}{p+1} \int_{ U_{-1}^+\cap U_{0}^+} |\ud(x)|^{p+1}dx + \frac{p-1}{p+1} \int_{U_0^+\cap U_1^+} |\ud(x)|^{p+1}dx \right) - \frac{1}{2\dx}\int_\R |\ud(x)|^{p+1}dx \\
	& \quad + \frac{1}{2\dx} \left(\frac{2}{p+1} \int_{U_{0}^-\cap U_{1}^-} |\ud(x)|^{p+1}dx + \frac{p-1}{p+1} \int_{U_{-1}^-\cap U_0^-} |\ud(x)|^{p+1}dx\right) \\
	& \le \frac{1}{2\dx} \int_{U_0^+} |\ud(x)|^{p+1}dx - \frac{1}{2\dx}\int_{\R} |\ud(x)|^{p+1}dx + \frac{1}{2\dx} \int_{U_0^-} |\ud(x)|^{p+1}dx
\end{align*}
and, hence, $I_1\le0$.

Thus, from \eqref{eq:ineq1112} we deduce:
\begin{align}
	\frac1p \frac{d}{dt} \|\ud(t)\|_p^p &\le \int_\R d^-_{\dx}\left(d^+_{\dx}\ud(x)\right) |\ud(x)|^{p-2}\ud(x)dx \label{eq:ineq11111} \\
	 &= -\frac{1}{\dx^2}\int_\R \left(u_\Delta(x+\dx)-u_\Delta(x)\right) \left(|\ud(x+\dx)|^{p-2}\ud(x+\dx) - |\ud(x)|^{p-2}\ud(x)\right)dx. \notag
\end{align}
Moreover, the following inequality (see \cite[Lemma II.5.5, p.~22]{Varopoulos:1992})
\begin{equation*}
	\left||x|^{p/2}-|y|^{p/2}\right|^2\le\frac{p^2}{4(p-1)}(x-y)\left(|x|^{p-2}x-|y|^{p-2}y\right),\quad\forall x,y\in\R,\quad \forall \ p \in (1,\infty),
\end{equation*}
guarantees that
\begin{equation}\label{eq:ineq1111}
	\frac{d}{dt} \|\ud(t)\|_p^p \le -\frac{4(p-1)}{p}\int_\R \left|\frac{|\ud(x+\dx)|^{p/2}-|\ud(x)|^{p/2}}{\dx}\right|^2=-\frac{4(p-1)}{p}\|d^+_{\dx}\big(|\ud|^{p/2}\big)\|_2^2 \le 0.
\end{equation}
This estimate and Lemma \ref{lemma:int1} allow us to write
\begin{equation}\label{eq:est22}
	\frac{d}{dt} \|\ud(t)\|_p^p +\frac{(p-1)}{p} \frac{\|\ud(t)\|_p^{p(p+1)/(p-1)}}{\|\ud^0(t)\|_1^{2p/(p-1)}} \le 0.
\end{equation}
Following the same arguments as in \cite{Escobedo:1991}, we conclude that for any $p\in[2,\infty)$
\begin{equation}\label{eq:est222}
	\|\ud(t)\|_p \le C(p)\|\ud^0\|_1\, t^{-\frac12(1-\frac1p)},\quad\forall t>0.
\end{equation}
The case $p\in(1,2)$ follows by interpolation. The case $p=\infty$ follows by tracking carefully the constants in \eqref{eq:est22} as in \cite{Veron:1979}.

Finally, the general case $\mu>0$ is immediate from \eqref{eq:est222} and the definition of $\umu$ \eqref{eq:umudef}, since for any $p\in[1,\infty]$ we have
\begin{equation*}
	\|\umu(t)\|_p = \mu^{1-\frac1p} \|\ud(\mu^2 t)\|_p \le C(p) \|\ud^0\|_1\, t^{-\frac12(1-\frac1p)}.
\end{equation*}
The proof is now complete.
\end{proof}

Now that we have estimates on the $L^p$-norms of the solution, we need to obtain a similar result for the discrete gradient.
\begin{prpstn}\label{prop:est2}
For all $p\in[1,\infty]$ there exists a constant $C=C(p,\|\ud^0\|_1 )>0$ such that:
\begin{equation}\label{eq:numgraddecay}
	\|d^+_{\dx/\mu} \umu(t)\|_p \le C t^{-\frac12(1-\frac1p)-\frac12},\quad \forall t>0,
\end{equation}
for all solutions of \eqref{eq:umu} with initial data $\ud^0\in L^1(\R)$.
\end{prpstn}
\begin{proof}
We proceed as in Proposition \ref{prop:graddecay}. Let us denote by $D^t_\mu$ the semigroup associated with
\begin{equation}\label{eq:linearkernelnum}
	\begin{cases}
		v_t(t,x)= \mu^2 \displaystyle\sum_{m=1}^{N} \omega_m v(t,x-m\frac{\dx}{\mu}) -\mu^2 \coefzero v(t,x) + \mu\coefone d^+_{\dx/\mu} v(t,x),& t>0,\mbox{ a.e. }x\in\R, \\[15pt]
		v(0,x)=v_0(x), & \mbox{ a.e. }x\in\R.
	\end{cases}
\end{equation}
Multiplying \eqref{eq:linearkernelnum} by $\sign(v(t,x))$, integrating on $\R$ and using that
\begin{equation*}
	\int_\R v(x-h)\sign(v(x)) dx \le \int_\R |v(x)| dx,\quad\forall h\in\R, 
\end{equation*}
one shows that $D^t_\mu$ is non-expansive in $L^1(\R)$.

Now, for every $\tau>0$ and $\mu>0$, the solution of \eqref{eq:umu} satisfies:
\begin{equation}\label{eq:intsemid}
	\umu(t+\tau)=G_\Delta^\mu(t)*D^t_\mu\umu(\tau)+\int_0^t G_\Delta^\mu(t-s)*D^{t-s}_\mu \Big(H\big(\umu(s+\tau)\big)\Big)ds,
\end{equation}
where
\begin{equation*}
	H(\umu(s,x))=\frac14\left(d^+_{\dx/\mu}\left(\umu(s,x)^2\right)+d^-_{\dx/\mu}\left(\umu(s,x)^2\right)\right) +\dx\,d^+_{\dx/\mu}R\Big(\umu(s,x-\frac{\dx}{\mu}),\umu(s,x)\Big)
\end{equation*}
and $G_\Delta^\mu$ is the fundamental solution of the one-dimensional semi-discrete heat equation, defined by
\begin{equation*}
	\big(G_\Delta^\mu(t)\big)_j=\frac{1}{2\pi}\int_{-\pi\mu/\dx}^{\pi\mu/\dx}e^{-\frac{4t\mu^2}{\dx^2}\sin^2\frac{\xi\dx}{2\mu}}e^{ij\xi\frac{\dx}{\mu}}d\xi,\quad j\in\Z,
\end{equation*}
It is well known (e.g. \cite{Brenner:1975}) that, for any $p\in[1,\infty]$,
\begin{equation*}
	\|G_\Delta^\mu(t)\|_p\le C(p)t^{-\frac12(1-\frac1p)},\quad t>0,
\end{equation*}
and
\begin{equation*}
	\|d^+_{\dx/\mu}G_\Delta^\mu(t)\|_p\le C(p)t^{-\frac12(1-\frac1p)-\frac12},\quad t>0.
\end{equation*}

Now let us apply the discrete operator $d^+_{\dx/\mu}$ to \eqref{eq:intsemid}. Then
\begin{equation}\label{eq:d+umu}
	d^+_{\dx/\mu}\umu(t+\tau)=d^+_{\dx/\mu}G_\Delta^\mu(t)*D^t_\mu\umu(\tau)+\int_0^t d^+_{\dx/\mu}G_\Delta^\mu(t-s)*D^{t-s}_\mu \Big(H\big(\umu(s+\tau)\big)\Big)ds.
\end{equation}
Using the decay properties of $G_\Delta^\mu$, Proposition \ref{prop:est1} and the $L^1$-stability of $D^t_\mu$, we obtain
\begin{align}\label{eq:estimd+}
	\|d^+_{\dx/\mu}\umu(t+\tau)\|_1 &\le \left\|d^+_{\dx/\mu}G_\Delta^\mu(t)\right\|_1 \left\|D^t_\mu\umu(\tau)\right\|_1 \\
	&\quad+\int_0^t \left\|d^+_{\dx/\mu}G_\Delta^\mu(t-s)\right\|_1\left\|D^{t-s}_\mu \Big(H\big(\umu(s+\tau)\big)\Big)\right\|_1 ds \notag \\
	&\le \left\|d^+_{\dx/\mu}G_\Delta^\mu(t)\right\|_1 \left\|\umu(\tau)\right\|_1+ \int_0^t \left\|d^+_{\dx/\mu}G_\Delta^\mu(t-s)\right\|_1\left\|H\big(\umu(s+\tau)\big)\right\|_1 ds \notag \\
	&\le C \left\|u_0\right\|_1 t^{-\frac12} + C \int_0^t (t-s)^{-\frac12} \left\|H\big(\umu(s+\tau)\big)\right\|_1 ds \notag.
\end{align}

We now prove that for any $p\in[1,\infty)$, we have
\begin{equation}\label{eq:estH}
	\|H(\umu(s+\tau))\|_p\le C\tau^{-\frac12}\|d^+_{\dx/\mu}\umu(s+\tau)\|_p
\end{equation}
Observe that, in view of Proposition \ref{prop:est1}, we have
\begin{equation*}
	\|u^\mu(s+\tau)\|_\infty\lesssim (s+\tau)^{-1/2}\|\umu(0)\|_1\lesssim \tau^{-1/2} \|u_0\|_1.
\end{equation*}
Thus we obtain
\begin{equation*}
	\left\| d^+_{\dx/\mu}\left(\umu(s+\tau)^2\right) \right\|_p \le2\left\|\umu(s+\tau)\right\|_\infty \left\|d^+_{\dx/\mu}\umu(s+\tau)\right\|_p \le C \tau^{-\frac12}\left\|d^+_{\dx/\mu}\umu(s+\tau)\right\|_p.
\end{equation*}
A similar result holds for $d^-_{\dx/\mu}$. Moreover, from the definition of $R$ in \eqref{eq:numviscR} we have:
\begin{align*}
	\dx \Big\| d^+_{\dx/\mu} &R\big(\umu(s+\tau,x-\frac{\dx}{\mu}),\umu(s+\tau,x)\big)\Big\|_p \le \frac{1}{2} \left\| d^+_{\dx/\mu}\left(\umu(s+\tau)|\umu(s+\tau)|\right)\right\|_p \\
	&\le \left\|\umu(s+\tau)\right\|_\infty \left\| d^+_{\dx/\mu} \umu(s+\tau) \right\|_p \le C \tau^{-\frac12} \left\| d^+_{\dx/\mu} \umu(s+\tau) \right\|_p,
\end{align*}
where we have used Proposition \ref{prop:est1} and that
\begin{equation*}
	\big|x|x|-y|y|\big|\le 2|x-y|\max\{|x|,|y|\},\quad\forall x,y\in\R.
\end{equation*}

Therefore, introducing \eqref{eq:estH} with $p=1$ in \eqref{eq:estimd+} we get
\begin{align*}
	\|d^+_{\dx/\mu}\umu(t+\tau)\|_1 &\le C t^{-\frac12} + C \tau^{-\frac12} \int_0^t (t-s)^{-\frac12} \left\| d^+_{\dx/\mu} \umu(s+\tau) \right\|_1 ds .
\end{align*} 
Applying fractional Gronwall's Lemma (see for example \cite[Lemma 2.4]{Bouharguane:2013}) and taking $t=\tau$, we conclude that
\begin{equation}\label{eq:estd+1}
	\|d^+_{\dx/\mu}(\umu(2\tau))\|_1 \le C_\tau,\quad \forall \mu>0,
\end{equation}
for some $C_\tau>0$ depending only on $\tau$ and $\|u_0\|_1$. It is enough now to use the definition of $\umu$ in \eqref{eq:umudef}, taking $\tau=1/2$ to obtain $\|d^+_{\dx}(\ud(\mu^2))\|_1 \le {C}/{\mu}$, for all $\mu>0$. Putting $\mu^2=t$ we find 
\begin{equation*}
	\|d^+_{\dx}(\ud(t))\|_1 \le \frac{C}{\sqrt{t}},\quad \forall t>0,
\end{equation*}
that is, \eqref{eq:numgraddecay} for $\mu=1$ and $p=1$.

The case $\mu=1$ and $p\in(1,\infty)$ is immediate from \eqref{eq:d+umu}, \eqref{eq:estH} and \eqref{eq:estd+1}. Indeed, we have
\begin{align*}
	\|d^+_{\dx/\mu}\umu(t+\tau)\|_p &\le \left\|d^+_{\dx/\mu}G_\Delta^\mu(t)\right\|_p \left\|\umu(\tau)\right\|_1+ \int_0^t \left\|d^+_{\dx/\mu}G_\Delta^\mu(t-s)\right\|_p\left\|H\big(\umu(s+\tau)\big)\right\|_1 ds \notag \\
	&\le C_\tau t^{-\frac12(1-\frac1p)-\frac12} + C_\tau \int_0^t (t-s)^{-\frac12(1-\frac1p)-\frac12} ds \notag.
\end{align*}
with $C_\tau=C(p,\tau,\|u_0\|_1)$. Applying again fractional Gronwall's Lemma and taking $t=\tau$ we obtain that
\begin{equation}\label{eq:estd+p}
	\|d^+_{\dx/\mu}\umu(2\tau)\|_p \le C_\tau,\quad \forall \mu>0,
\end{equation}
This is equivalent to \eqref{eq:numgraddecay} for $\mu=1$ and $p\in(1,\infty)$.

Furthermore, repeating similar arguments, the case $p=\infty$ follows from \eqref{eq:d+umu} and estimates \eqref{eq:estH} and \eqref{eq:estd+p}:
\begin{align*}
	\|d^+_{\dx/\mu}\umu(t+\tau)\|_\infty &\le \left\|d^+_{\dx/\mu}G_\Delta^\mu(t)\right\|_\infty \left\|\umu(\tau)\right\|_1+ \int_0^t \left\|d^+_{\dx/\mu}G_\Delta^\mu(t-s)\right\|_2\left\|H\big(\umu(s+\tau)\big)\right\|_{2} ds \\
	&\le C_\tau t^{-1} + C_\tau \int_0^t (t-s)^{-\frac34}ds.
\end{align*}
where $C_\tau=C(\tau,\|u_0\|_1)$. It is now enough to take $t=\tau$ to conclude that
\begin{equation*}
	\|d^+_{\dx/\mu}\umu(2\tau)\|_\infty \le C_\tau,\quad \forall \mu>0,
\end{equation*}
which is equivalent to \eqref{eq:numgraddecay} for $\mu=1$ and $p=\infty$.

Finally, the general case $\mu>0$ is immediate from the case $\mu=1$ and the definition of $\umu$ \eqref{eq:umudef}, since for any $p\in[1,\infty]$ we have
\begin{equation*}
	\|d^+_{\dx/\mu}\umu(t)\|_p = \mu^{2-\frac1p} \|d^+_{\dx}\ud(\mu^2 t)\|_p \le C t^{-\frac12(1-\frac1p)-\frac12}
\end{equation*}
This concludes the proof.
\end{proof}

%%%%%%%%%%% SUBSECTION 3.2 %%%%%%%%

\subsection{Compactness of the set $\{\umu\}_{\mu>0}$}
In this section, we prove the compactness of the trajectories of the family $\{\umu(t)\}_{\mu>0}$ introduced in the previous section, in order to pass to the limit $\mu\to\infty$. Unlike the continuous case, we do not have estimates of $\umu$ in $H^1(\R)$, since $\umu$ is piecewise constant. Nevertheless, the following lemma makes possible the use of the compact embedding of $H^s_{loc}(\R)$ into $L^2_{loc}(\R)$, with $0<s<1/2$. The proof will be given in the Appendix. 
\begin{lmm}\label{lemma:hs}
	For any $0<s<\frac12$, there exists a positive constant $C=C(s)$ such that, for any mesh-size $0<\dx<1$, the following holds for all piecewise constant functions $w$ as in \eqref{eq:udelta}:
	\begin{equation*}
		\|w\|_{H^s(\R)}\le C \left( \|w\|_2+\|d^+_{\dx}w\|_2\right).
	\end{equation*}
\end{lmm}
Let us remark that, as a consequence of this lemma and Proposition \ref{prop:est1} and Proposition \ref{prop:est2}, we obtain a time-decay estimate for $\umu$ in $H^s(\R)$ with $0<s<1/2$ (this can be done since $\umu$ is piecewise constant on intervals of length $\dx/\lambda$):
\begin{equation}\label{eq:hsbound}
	\|\umu(t)\|_{H^s(\R)}\le C \left( \|\umu(t)\|_2+\|d^+_{\dx/\mu}\umu(t)\|_2\right)\le C \left(t^{-\frac14}+t^{-\frac34} \right),\quad\forall t>0,\forall\mu>0.
\end{equation}
Thus, we can use Theorem \ref{teo:simon} to prove the compactness of the family $\{\umu\}_{\mu>0}$.
\begin{thrm}\label{teo:comp}
	For every $0<\tau<T<\infty$, the family $\{\umu\}_{\mu>0}\subset C([\tau,T],L^1(\R))$ is relatively compact.
\end{thrm}
\begin{proof} We will proceed in two steps, analogously to Theorem \ref{teo:compact}.

\textit{Step 1.} First we will show the result locally in $C([\tau,T],L^1(I))$, with $I=[-r,r]$ for an arbitrary $r>0$. Let us consider the spaces $X=H^s(I)$ with $s\in(0,\frac12)$, $Z=L^2(I)$ and $Y=H^{-1}(I)$.

From \eqref{eq:hsbound} we know that the set $\{\umu\}_{\mu>0}$ is bounded in $L^\infty([\tau,T],H^s_{loc}(\R))$. In particular, the first condition of Theorem \ref{teo:simon} is fulfilled. Thus, it suffices to check that $\umu_t$ is bounded in $L^\infty([\tau,T],H^{-1}(I))$. Let us multiply \eqref{eq:umu} by a function $\varphi\in C^\infty_c(\R)$ and integrate it over $\R$. Using the definition of $R$ in \eqref{eq:numviscR}, we have:
\begin{align*}
	\left|\int_\R \umu_t \varphi dx\right| &\le \frac{1}{4}\left|\int_\R \left(d^+_{\dx/\mu} \big(\umu(x)^2\big)+d^-_{\dx/\mu} \big(\umu(x)^2\big)\right)\varphi(x) dx\right| \\
	&\quad + \dx \left|\int_\R d^+_{\dx/\mu} R\big(\umu(x-\frac{\dx}{\mu}),\umu(x)\big) \varphi(x) dx\right| + \left|\int_\R d^-_{\dx/\mu}\big(d^+_{\dx/\mu} \umu(x)\big)\varphi(x) dx\right| \\
	&\quad +\left|\int_\R\left(\mu^2\sum_{m=1}^{N} \omega_m \umu(x-m\frac{\dx}{\mu}) -\mu^2 \coefzero \umu(x) +\mu \coefone d^+_{\dx/\mu} \umu(x)\right)\varphi(x) dx\right| \\
	&\le \frac{1}{2}\|d^+_{\dx} \varphi\|_2 \|\umu\|_4^2 + \frac{1}{2}\|d^+_{\dx} \varphi\|_2 \|\umu\|_4^2 +\|d^+_{\dx} \varphi \|_2 \|d^+_{\dx} \umu \|_2 \label{eq:comp1}\\
	&\quad +\left|\int_\R\left(\mu^2 \sum_{m=1}^{N} \omega_m \big(\umu(x-m\frac{\dx}{\mu})- \umu(x) \big) + \mu \coefone d^+_{\dx/\mu}\umu(x)\right)\varphi(x) dx\right|.
\end{align*}
Obviously, the first three terms on the right hand side of the inequality are uniformly bounded for $\mu>0$, so let us focus on the last one. Using the Fourier transform and the definition of $\coefzero$ in \eqref{eq:factorfn}, we have
\begin{align*}
	I_{\mu}&=\left|\int_\R\left(\mu^2 \sum_{m=1}^{N} \omega_m \big(\umu(x-m\frac{\dx}{\mu})- \umu(x) \big) + \mu \coefone d^+_{\dx/\mu}\umu(x)\right)\varphi(x) dx\right| \\
	&\le \mu^2 \int_\R \left| \sum_{m=1}^{N} \omega_m \left(e^{-im\frac{\dx}{\mu}\xi}- 1\right)+ \coefone \frac{e^{i\frac{\dx}{\mu}\xi}-1}{\dx}\right| \left|\widehat{\umu}(\xi)\right|\left|\widehat{\varphi}(\xi)\right| d\xi.
\end{align*}
If we take $a=e^{-\dx}$ and $b=e^{-i\frac{\dx}{\mu}\xi}$ on Lemma \ref{lemma:serie} and use the definitions of $\omega_m$ in \eqref{eq:factorwm} and $\coefone$ in \eqref{eq:factorfn}, we have:
\begin{align}
	\Bigg|\sum_{m=1}^{N} \omega_m & \left(e^{-im\frac{\dx}{\mu}\xi}- 1\right) + \coefone \frac{e^{i\frac{\dx}{\mu}\xi}-1}{\dx} \Bigg| \\
	& = \left|e^{\dx}-1\right| \left| \sum_{m=1}^{N} e^{-m\dx} \left(e^{-im\frac{\dx}{\mu}\xi}- 1\right) + \sum_{m=1}^{N} m e^{-m\dx} (e^{i\frac{\dx}{\mu}\xi}-1)\right| \notag \\
	& \le \left|e^{\dx}-1\right| \left|e^{-i\frac{\dx}{\mu}\xi}-1\right|^2 \frac{e^{-\dx}}{(1-e^{-\dx})^3}  =   \frac{|e^{-i\frac{\dx}{\mu}\xi}-1|^2}{(1-e^{-\dx})^2}. \notag
\end{align}
Combining this inequality with the Cauchy-Schwartz inequality and the fact that
\begin{equation*}
	\|d^+_{\dx/\mu}\umu\|_2^2=\int _{\R}\Big|\frac{e^{i\frac{\dx}{\mu}\xi}-1}{\dx/\mu} \Big|^2 |\widehat{\umu}(\xi)|^2d\xi.
\end{equation*}
we obtain
\begin{equation*}
	I_{\mu} \le \frac{\dx^2}{(1-e^{-\dx})^2} \|d^+_{\dx/\mu} \umu\|_2 \|d^+_{\dx/\mu} \varphi\|_2.
\end{equation*}
Thus, using that $\|d^+_{\dx/\mu}\varphi\|_2 \lesssim \|\varphi'\|_2$, we get
\begin{align*}
	\left|\int_\R \umu_t(t)\varphi dx\right|&\le \|d^+_{\dx/\mu} \varphi \|_2 \|\umu(t)\|_4^2 +\|d^+_{\dx/\mu} \varphi \|_2 \|d^+_{\dx/\mu} \umu(t) \|_2 +\frac{\dx^2}{(1-e^{-\dx})^2} \|d^+_{\dx/\mu} \umu(t)\|_2 \|d^+_{\dx/\mu} \varphi\|_2 \\
	&\le C \|\varphi \|_{H^1(\R)} \Big(\|\umu(t)\|_4^2 +\|d^+_{\dx/\mu} \umu(t) \|_2\Big).
\end{align*}
for any $\varphi\in C^\infty_c(I)$ and with $C>0$ independent of $\mu$. In view of Propositions \ref{prop:est1} and \ref{prop:est2}, both norms of $\umu$ in the right-hand side are uniformly bounded in $[\tau, T]$, so $\umu_t$ is uniformly bounded in $L^\infty([\tau,T],H^{-1}(I))$. We conclude that the family $\{\umu\}_{\mu>0}$ is relatively compact in $C([\tau,T],L^2_{loc}(\R))$. Finally, compactness in $L^2_{loc}(\R)$ implies compactness in $L^1_{loc}(\R)$, so $\{\umu\}_{\mu>0}$ is also relatively compact in $C([\tau,T],L^1_{loc}(\R))$.

\textit{Step 2.} Now we need to extend the result globally. Let us consider again the same function $\psi_r$ defined in the third step of the proof of Theorem \ref{teo:compact}, such that $\psi_r(z)=\psi(z/r)$ with $\psi$ given by \eqref{eq:psi} and $r>0$. Since we know that $\{\umu\}_{\mu>0}$ is relatively compact in $C([\tau,T],L^1_{loc}(\R))$, it suffices to show that
\begin{equation}\label{eq:sup2}
	\sup_{[\tau,T]}\|\umu(t)\psi_r\|_1\longrightarrow0 \quad\mbox{as }r\rightarrow\infty,\mbox{ uniformly on } \mu\ge1.
\end{equation}
A similar argument as in Theorem \ref{teo:compact} shows that it is enough to prove \eqref{eq:sup2} for nonnegative initial data and solutions. Thus, we focus on those. Let us multiply \eqref{eq:umu} by $\psi_r$ and integrate it over $(0,t)\times\R$. We obtain:
\begin{align}
	\int_\R &\umu(t,x)\psi_r(x) dx = \int_\R \umu_0(x)\psi_r(x) dx \\
	& +\frac{1}{4}\int_0^t\int_\R \left(d^+_{\dx/\mu}\big(\umu(s,x)^2\big) + d^-_{\dx/\mu}\big(\umu(s,x)^2\big) \right)\psi_r(x)dxds \notag \\
	& + \dx\int_0^t\int_\R d^+_{\dx/\mu} \Big(R\big(\umu(s,x-\frac{\dx}{\mu}),\umu(s,x)\big)\Big)\psi_r(x)dxds \notag \\
	& +\int_0^t\int_\R d^-_{\dx/\mu} \Big(d^+_{\dx/\mu} \big(\umu(s,x)\big)\Big)\psi_r(x)dxds \notag \\
	& +\int_0^t \int_\R \left(\mu^2\sum_{m=1}^{N} \omega_m \big(\umu(s,x-m\frac{\dx}{\mu}) - \umu(s,x)\big) +\mu \coefone d^+_{\dx/\mu} \umu(s,x)\right) \psi_r(x)dxds. \notag
\end{align}
We pass now the discrete derivatives to $\psi_r$ and estimate the right-hand side using time-decay estimates from Proposition \ref{prop:est1}:
\begin{align}
	&\int_\R \umu(t,x)\psi_r(x) dx \lesssim \int_\R \umu_0(x)\psi_r(x) dx +\|\psi'\|_\infty \frac{\sqrt{t}}{r} +\|\psi''\|_\infty \frac{t}{r^2} \label{eq:comp2}\\
	&\quad +\int_0^t \int_\R \left(\mu^2\sum_{m=1}^{N} \omega_m \big(\umu(s,x-m\frac{\dx}{\mu}) - \umu(s,x)\big) +\mu \coefone d^+_{\dx/\mu} \umu(s,x)\right) \psi_r(x)dxds, \notag
\end{align}
Let us focus on the last term. Using Taylor expansions and the definition of $\coefzero$ and $\coefone$ from \eqref{eq:factorfn}, we have
\begin{align*}
	\int_\R \Big(\mu^2\sum_{m=1}^{N} &\omega_m \big(\umu(s,x-m\frac{\dx}{\mu}) - \umu(s,x)\big) +\mu \coefone d^+_{\dx/\mu} \umu(s,x)\Big) \psi_r(x)dx \\ 
	&= \mu^2 \sum_{m=1}^{N} \omega_m \int_\R \umu(s,x) \Big(\psi_r(x+m\frac{\dx}{\mu}) - \psi_r(x)- m\frac{\dx}{\mu}\, d^-_{\dx/\mu} (\psi_r(x))\Big) dx \\
	&\lesssim \left\|\mu^2\sum_{m=1}^{N} \omega_m \Big(\psi_r(x+m\frac{\dx}{\mu}) - \psi_r(x)- m\frac{\dx}{\mu}\, d^-_{\dx/\mu} (\psi_r(x))\Big)\right\|_\infty \|\umu(s)\|_1 \\
	&\lesssim \frac{\|\psi''\|_\infty}{r^2} \|\ud^0\|_1.
\end{align*}
Thus, plugging this into \eqref{eq:comp2} and using the non-negativity of the solution, we get
\begin{align}
	\int_\R |\umu(t,x)|\psi_r(x) dx \lesssim \int_\R |\umu_0(x)|\psi_r(x) dx +\|\psi'\|_\infty \frac{\sqrt{t}}{r} +\|\psi''\|_\infty \frac{t}{r^2},
\end{align}
which tends to $0$ uniformly on $\mu>0$ when $r\rightarrow\infty$. Therefore, we proved \eqref{eq:sup2} and, consequently, we can assure that $\{\umu\}_{\mu>0}$ is relatively compact in $C([\tau,T],L^1(\R))$.
\end{proof}

A slight modification of the proof of the previous theorem gives as the necessary estimate to identify the initial data, stated in the following proposition.
\begin{lmm}\label{lemma:compini}
For every test function $\varphi\in C^\infty_c(\R)$, there exists $C>0$, independent of $\mu$, such that
	\begin{equation}\label{eq:compini}
		\left|\int_\R \umu(t,x)\varphi(x)dx-\int_\R \umu_0(x)\varphi(x)dx\right|\le C(t+\sqrt{t}).
	\end{equation}
\end{lmm}
\begin{proof}
It is enough to multiply \eqref{eq:umu} by $\varphi\in C^\infty_c(\R)$ and integrate it over $(0,t)\times\R$. Then, integrating by parts and repeating arguments similar to the ones in the second step of the proof for Theorem \ref{teo:comp}, we deduce \eqref{eq:compini}.
\end{proof}

%%%%%%%%%%% SUBSECTION 3.3 %%%%%%%%

\subsection{Passing to the limit}
Finally, we have everything that we need to prove our main result, stated in Theorem \ref{teo:asymptoticsnum}, regarding the large-time behavior of the approximations to the solution of problem \eqref{eq:model}.
\begin{proof}[Proof of Theorem \ref{teo:asymptoticsnum}]
By Theorem \ref{teo:comp}, we know that for every $0<\tau<T<\infty$, the family $\{\umu\}_{\mu>0}$ is relatively compact in $C([\tau,T],L^1(\R))$. Consequently, there exists a subsequence of it (which we will not relabel) and a function $\bar u\in C((0,\infty),L^1(\R))$ such that
\begin{equation}\label{eq:sequence2}
	\umu\longrightarrow \bar u\in C([\tau,T],L^1(\R)),\quad\mbox{as }\mu\rightarrow\infty.
\end{equation}
We can also assume that $\umu(t,x)\rightarrow\bar u(t,x)$ almost everywhere in $(0,\infty)\times\R$ as ${\mu\rightarrow\infty}$.

Now, we multiply equation \eqref{eq:umu} by a test function $\phi\in C^\infty_c((0,\infty)\times\R)$ and integrate it over $(0,\infty)\times\R$. We have:
\begin{align}\label{eq:umuweak}
		\int_0^\infty \int_\R &\umu_t (t,x)\phi(t,x)dxdt=\frac{1}{4}\int_0^\infty \int_\R \Big(d^+_{\dx/\mu} \big(\umu(t,x)^2\big)+d^-_{\dx/\mu} \big(\umu(t,x)^2\big)\Big) \phi(t,x) dx dt \\
		&\quad + \dx \int_0^\infty \int_\R d^+_{\dx/\mu} R\big(\umu(t,x-\frac{\dx}{\mu}),\umu(t,x)\big) \phi(t,x) dx dt \notag \\
		&\quad + \int_0^\infty \int_\R d^-_{\dx/\mu} \left(d^+_{\dx/\mu} \umu(t,x) \right) \phi(t,x) dx dt \notag \\
		&\quad + \int_0^\infty \int_\R \left( \mu^2\sum_{m=1}^{N} \omega_m \umu(t,x-m\frac{\dx}{\mu}) -\mu^2 \coefzero \umu(t,x) +\mu \coefone d^+_{\dx/\mu} \umu(t,x) \right) \phi(t,x) dx dt \notag
\end{align}
Our claim is that, passing to the limit $\mu\rightarrow\infty$, we obtain that $\bar u$ is a weak solution of the equation:
\begin{equation}\label{eq:burgersvisc2}
	\begin{cases}
		\bar u_{t}=\bar u \bar u_{x}+(1+\coeftwo )\bar u_{xx},&(t,x)\in(0,\infty)\times\R,\\
		\bar u(0)=M\delta_0.
	\end{cases}
\end{equation}
All the limits in \eqref{eq:umuweak} are known except for the last term. In fact, let us recall that the degree of homogeneity of $R$ makes the corresponding numerical viscosity term vanish as $\mu\to\infty$, as shown in \cite{Pozo:2015}.

Thus, it is sufficient to check that we can take the limit $\mu\to\infty$ in
\begin{equation*}
	\mathcal L^\mu(t)= \int_\R\left(\mu^2\sum_{m=1}^{N} \omega_m \umu(t,x-m\frac{\dx}{\mu}) -\mu^2 \coefzero \umu(t,x) +\mu \coefone d^+_{\dx/\mu} \umu(t,x)\right)\phi(t,x) dx.
\end{equation*}
First, we reorder $\mathcal{L}_\mu$:
\begin{equation} \label{eq:convlap2}
	\mathcal L^\mu(t)=\mu^2\int_\R \umu(t,x) \sum_{m=1}^N \omega_m \left(\phi(t,x+m\frac{\dx}{\mu})-\phi(t,x)-m\frac{\dx}{\mu} d^-_{\dx/\mu}(\phi(t,x))\right) dx.
\end{equation}
Now, due to Taylor's Theorem, for each $m\in\{1,\dots,N\}$
\begin{equation*}
	\phi(t,x+m\frac{\dx}{\mu})-\phi(t,x)=m\frac{\dx}{\mu}\phi_x(t,x)+\frac{1}{2}m^2\frac{\dx^2}{\mu^2}\phi_{xx}(t,x)+\frac1{\mu^3} O(\|\phi_{xxx}(t)\|_\infty).
\end{equation*}
In the same way, 
\begin{equation*}
	d^-_{\dx/\mu}(\phi(t,x))=\phi_x(t,x)-\frac{1}{2}\frac{\dx}{\mu}\phi_{xx}(t,x)+\frac1{\mu^2}O(\|\phi_{xxx}(t)\|_\infty).
\end{equation*}
We combine this into \eqref{eq:convlap2} and get
\begin{equation} \label{eq:convfin2}
	 \mathcal L^\mu(t)= \coeftwo \int_\R \umu(t,x)\phi_{xx}(t,x) dx +O(\|\phi_{xxx}(t)\|_\infty) \frac{1}{\mu} \int_\R \umu(t,x) dx,
\end{equation}
where
\begin{equation*}
	\coeftwo =\frac{\dx^2}{2} \left(\sum_{m=1}^N m(m-1)\omega_m \right).
\end{equation*}
Therefore, as $\umu\rightarrow\bar u$ in $C([\tau,T],L^1(\R))$, taking the limit $\mu\to\infty$ in \eqref{eq:convfin2}, we obtain:
\begin{equation*}
	\lim_{\mu\rightarrow\infty} \int_0^\infty \mathcal L^\mu(t)= \coeftwo \int_0^\infty \int_\R \bar u(t,x)\phi_{xx}(t,x)dx.
\end{equation*}

\begin{rmrk}
	Let us emphasize the key role that the correcting factors $\coefzero$ and $\coefone$ have on the limit above. Note that the fact that
$$\int_\R K(z) dz = \int z K(z) dz = 1$$
played an important role in the proof of Theorem \ref{teo:compact}, allowing us to show that $K*u_{xx}$ behaves as $u_{xx}$ as $t\to\infty$. Moreover, this is related also with the decomposition of $K$ in Dirac delta functions as in \cite{Duoandikoetxea:1992}. Now, at the discrete level, the corrector factors $\coefzero $ and $\coefone $ had to be chosen accordingly. In this case, due to the truncation of the integral, the use of a third factor $\coeftwo$ is required, though. All the same, these three coefficients are the discretized moments of the kernel $K$. 
\end{rmrk}

It follows that $\bar u$ satisfies
\begin{equation*}
	-\int_0^\infty \int_\R \bar u \phi_t = -\frac12 \int_0^\infty \int_\R \bar u^2 \phi_x + (1+\coeftwo ) \int_0^\infty \int_\R \bar u \phi_{xx}, 
\end{equation*}
so it is a weak solution of the equation in \eqref{eq:burgersvisc2}. It remains to identify the behavior of $\bar u$ as $t\to0$. Due to Lemma \ref{lemma:compini}, for any $\varphi\in C^\infty_c(\R)$ we have
\begin{equation*}
	\left|\int_\R \umu(t,x)\varphi(x)dx - \int_\R \umu_0(x)\varphi(x)dx \right| \le C(t+\sqrt{t})
\end{equation*}
and from \eqref{eq:sequence2} we deduce
\begin{equation*}
	\left|\int_\R \bar u(t,x)\varphi(x)dx - M\varphi(0)\right| \le C(t+\sqrt{t})
\end{equation*}
by letting $\mu\to\infty$. Passing to the limit $t\to0$ and using classical approximation arguments, we deduce that $\bar u(0)=M\delta_0$ in the sense of bounded measures. Thus, we conclude that $\bar u$ is the unique solution $u_M$ of equation \eqref{eq:burgersvisc2}, and that, in fact, the whole family $\{\umu\}_{\mu>0}$ converges to $u_M$ in $C((0,\infty),L^1(\R))$.

Therefore, by \eqref{eq:sequence2}, we have:
\begin{equation*}
	\lim_{\mu\rightarrow\infty}\|\umu(1)-u_M(1)\|_1=0
\end{equation*}
and setting $\mu=\sqrt{t}$ and making use of the self-similar form of $u_M$ (see e.g. \cite{Escobedo:1991}) we obtain
\begin{equation}\label{eq:limnorm12}
	\lim_{t\rightarrow\infty}\|\ud(t)-u_M(t)\|_1=0.
\end{equation}

Finally, the convergence in the $L^p$-norms for $p\in(1,\infty)$ follows from \eqref{eq:limnorm12}, the decay estimate of Proposition \ref{prop:est1} for $p=\infty$ and the H\"older inequality. In fact, we have:
\begin{equation*}
	\|\ud(t)-u_M(t)\|_p\le\left(\|\ud(t)\|_\infty+\|u_M(t)\|_\infty\right)^{1-\frac1p}\|\ud(t)-u_M(t)\|_1^{\frac1p}\le o(t^{-\frac12(1-\frac1p)}).
\end{equation*}
Using the piecewise constant interpolation of $u_M$, which we denote $S(u_M)$, and \eqref{eq:lemmapcest} from the Appendix, the case $p=\infty$ follows:
\begin{align*}
	\|\ud(t)-u_M(t)\|_\infty &\le \|\ud(t)-S(u_M(t))\|_\infty + \|S(u_M(t))-u_M(t)\|_\infty \\
	& \lesssim \|\ud(t)-S(u_M(t))\|_2^{\frac12} \|d^+_{\dx} (\ud(t)-S(u_M(t)))\|_2^{\frac12} + \dx \|u_{M,x}(t)\|_\infty \\
	& \lesssim (\|\ud(t)-u_M(t)\|_2 + \|u_M(t)-S(u_M(t))\|_2)^{\frac12} (\|d^+_{\dx} \ud(t)\|_2 + \|d^+_{\dx} S(u_M(t))\|_2)^{\frac12} \\
	&\qquad + \dx \|u_{M,x}(t)\|_\infty \\
	& \le o(t^{-\frac12} + t^{-\frac34} + t^{-1}).
\end{align*}
Now the proof is complete.
\end{proof}

%%%%%%%%%%% SUBSECTION 3.4 %%%%%%%%

\subsection{Convergence of the scheme}
To conclude this section, let us prove that $u_\Delta$ converges to the solution $u$ of \eqref{eq:model} as $\dx\to0$.

\begin{thrm}\label{teo:convergence}
	Let $u_0\in L^1(\R)$ and $N=N(\dx)\in\N$ such that $N\dx\to\infty$ as $\dx\to0$. The set of approximated solutions $\{\ud\}_{\dx>0}$ given by \eqref{eq:scheme} converges in $C((0,\infty),L^1(\R))$ to the solution $u$ of \eqref{eq:model} as $\dx\to0$.
\end{thrm}
\begin{proof}
Following the same arguments as in Theorem \ref{teo:comp}, one shows that for every $0<\tau<T<\infty$, the family $\{\ud\}_{\dx>0}\subset C([\tau,T],L^1(\R))$ is relatively compact. Thus, there exists a subsequence of it (which we will not relabel) and a function $\bar u\in C((0,\infty),L^1(\R))$ such that
\begin{equation}\label{eq:sequence3}
	\ud\longrightarrow \bar u\in C([\tau,T],L^1(\R)),\quad\mbox{as }\dx\rightarrow0.
\end{equation}
We can also assume that $\ud(t,x)\rightarrow\bar u(t,x)$ almost everywhere in $(0,\infty)\times\R$ as ${\dx\rightarrow0}$.

Now, we take $\mu=1$ in equation \eqref{eq:umu}, multiply it by a test function $\phi\in C^\infty_c((0,\infty)\times\R)$ and integrate it over $(0,\infty)\times\R$. We have:
\begin{align}\label{eq:umuweak2}
		\int_0^\infty \int_\R &u_{\Delta,t} (t,x)\phi(t,x)dxdt=\frac{1}{4}\int_0^\infty \int_\R \Big(d^+_{\dx} \big(\ud(t,x)^2\big)+d^-_{\dx} \big(\ud(t,x)^2\big)\Big) \phi(t,x) dx dt \\
		&\quad + \dx \int_0^\infty \int_\R d^+_{\dx} R\big(\ud(t,x-\dx),\ud(t,x)\big) \phi(t,x) dx dt \notag \\
		&\quad + \int_0^\infty \int_\R d^-_{\dx} \left(d^+_{\dx} \ud(t,x) \right) \phi(t,x) dx dt \notag \\
		&\quad + \int_0^\infty \int_\R \left( \sum_{m=1}^{N} \omega_m \ud(t,x-m\dx) - \coefzero \ud(t,x) + \coefone d^+_{\dx} \ud(t,x) \right) \phi(t,x) dx dt \notag
\end{align}
Our claim is that, passing to the limit $\dx\rightarrow0$, we obtain that $\bar u$ is a weak solution of the equation \eqref{eq:model}. Using classical arguments, thanks to \eqref{eq:sequence3}, Proposition \ref{prop:est1} and Proposition \ref{prop:est2}, we can take all the limits in \eqref{eq:umuweak2}, except for the last term. Thus, it is sufficient to check that we can also pass to the limit $\dx\to0$ in
\begin{align*}
	\mathcal L_\Delta(t) & = \int_\R\left(\sum_{m=1}^{N} \omega_m \ud(t,x-m\dx) - \coefzero \ud(t,x) + \coefone d^+_{\dx} \ud(t,x)\right)\phi(t,x) dx \\
	& = \int_\R \ud(t,x) \left(\sum_{m=1}^{N} \omega_m \phi(t,x+m\dx) - \coefzero \phi(t,x) - \coefone d^-_{\dx} \phi(t,x)\right) dx.
\end{align*}
First, let us first observe that
\begin{equation*}
	\coefzero = \sum_{m=1}^N e^{-m\dx} (e^{\dx}-1) = 1-e^{-N\dx}\to1
\end{equation*}
and
\begin{equation*}
	\coefone =\dx (e^{\dx}-1) \sum_{m=1}^N m e^{-m\dx} = \frac{\dx e^{-N\dx}(e^{(N+1)\dx}-e^{\dx}(N+1)+N)}{e^{\dx}-1}\to1,
\end{equation*}
as long as $N=N(\dx)$ is taken such that $N\dx\to\infty$ as $\dx\to0$. Moreover, using \eqref{eq:fourierK} and that
\begin{equation*}
	(e^{\dx}-1) \frac{1-e^{-N\dx(1-i\xi)}}{e^{(1-i\xi)\dx}-1} \to \frac{1}{1-i\xi},\quad\mbox{as }\dx\to0
\end{equation*}
we obtain
\begin{align*}
	\Bigg|\sum_{m=1}^{N} \omega_m \phi(t,x+m\dx) - \tilde K*\phi(t,x) \Bigg| &\le \int_\R |\widehat{\phi}(t,\xi)| \left| \sum_{m=1}^{N} \omega_m e^{im\dx\xi} - \widehat{K}(-\xi)\right| d\xi \\
	& = \int_\R |\widehat{\phi}(t,\xi)|\left|(e^{\dx}-1) \frac{1-e^{-N\dx(1-i\xi)}}{e^{(1-i\xi)\dx}-1}-\frac{1}{1-i\xi}\right| d\xi \to 0,
\end{align*}
where $\tilde K(z)=K(-z)$. Therefore
\begin{equation*}
	\lim_{\dx\to0} \mathcal L_\Delta(t) = \int_\R \bar u(t,x) \left(\tilde K*\phi(t,x) - \phi(t,x) - \phi_x(t,x)\right) dx.
\end{equation*}
It follows that $\bar u$ satisfies
\begin{equation*}
	-\int_0^\infty \int_\R \bar u \phi_t = -\frac12 \int_0^\infty \int_\R \bar u^2 \phi_x + \int_0^\infty \int_\R \bar u \phi_{xx} + \int_0^\infty \int_\R \bar u (\tilde K*\phi-\phi-\phi_x).
\end{equation*}
so it is a weak solution of the equation in \eqref{eq:model}.

Now, it remains to identify the behavior of $\bar u$ as $t\to0$. Using similar estimates as in the proof for Lemma \ref{lemma:compini}, we can show that for every test function $\varphi\in C^\infty_c(\R)$ and $\dx<1$, there exists $C>0$, independent of $\dx$, such that
\begin{equation*}
	\left|\int_\R \ud(t,x)\varphi(x)dx-\int_\R \ud^0(x)\varphi(x)dx\right|\le C(t+\sqrt{t}).
\end{equation*}
and from \eqref{eq:sequence3} and the definition of $\ud^0$ in \eqref{eq:udelta}, we deduce
\begin{equation*}
	\left|\int_\R \bar u(t,x)\varphi(x)dx - \int_\R u_0(x)\varphi(x)dx\right| \le C(t+\sqrt{t})
\end{equation*}
by letting $\dx\to0$. Using an approximation argument we deduce that $u(t)\rightarrow u_0$ in $L^1(\R)$ as $t\rightarrow 0$. Thus, we conclude that $\bar u$ is the unique solution $u$ of equation \eqref{eq:model} and that, in fact, the whole family $\{\ud\}_{\dx>0}$ converges to $u$ in $C((0,\infty),L^1(\R))$. Now the proof is complete.
\end{proof}

%%%%%%%%%%%% SECTION 4 %%%%%%%%%%%%%%%%%%%%%%%%%

\section{Numerical experiments}\label{section4}
The aim of this last section is to support the necessity of using large-time behavior preserving schemes for the augmented Burgers equation. On the one hand, we show the importance of a numerical flux that does not destroy the N-wave shape at the early stages. On the other, we emphasize the role of the correcting factors $\coefzero$ and $\coefone$ in the truncation of the convolution. Note that the former phenomenon was already stated in \cite{Pozo:2015} in the hyperbolic case, while the latter is an original contribution of the present work.

Regarding the time discretization, we opt for the explicit Euler for its simplicity. Even if there is no guarantee that the asymptotic behavior is preserved, numerical simulations exhibit a correct performance. Thus, we consider it enough to illustrate the key points enumerated above. We need to take into account that there is a stability condition that must be satisfied to ensure the convergence. It is easy to see (e.g. \cite{Castro:2010,Godlewski:1991}) that a sufficient condition is that
\begin{equation}
	\frac{\dt}{\dx}\max_{j}\{u^0_j\}+2\nu\frac{\dt}{\dx^2} +c\,\dt \sum_{m=1}^{N} (m+1) \omega_m \le 1.
\end{equation}

Let us choose the following compactly supported initial data.
\begin{equation}
	u_0(x)=\begin{cases} -\dfrac{1}{10}\sin\left(\dfrac{x}{2}\right),&x\in[-\pi,0], \\[5pt] -\dfrac{1}{20}\sin(2x), &x\in[0,\dfrac{\pi}{2}], \\[5pt] 0,&\mbox{elsewhere}\end{cases}
\end{equation}
We take a mesh size $\dx=0.1$. In order to avoid boundary issues, we choose a large enough spatial domain.

\begin{figure}[ht!]
	\includegraphics[width=0.5\linewidth]{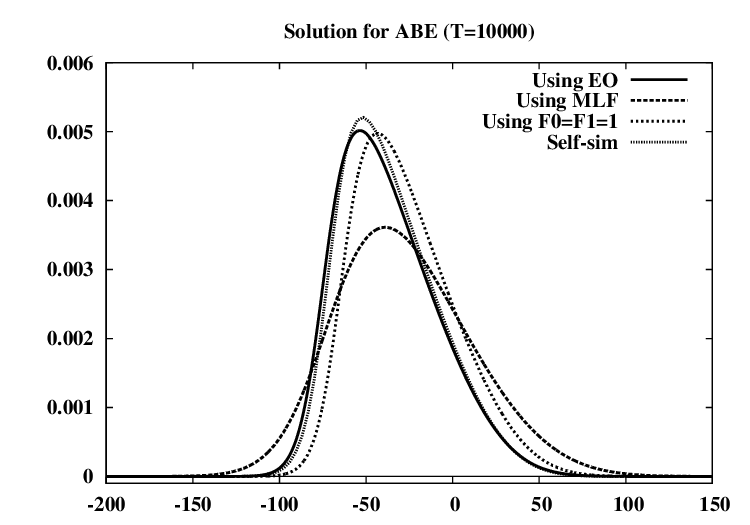}
	\caption{Solution of ABE with $\nu=10^{-2}$, $c=2\times10^{-2}$ and $\theta=1$ at $t=10^4$, using scheme \eqref{eq:scheme} discretized explicitly. We use EO (solid) and modified LF (dashed) numerical fluxes for the nonlinearity, as well as EO without correcting factors (dotted), comparing the solutions to the asymptotic profile (gray).}
	\label{fig1}
\end{figure}

In Figure \ref{fig1} we show the solution for $\nu=10^{-2}$, $c=2\times10^{-2}$ and $\theta=1$ at time $t=10^4$, as well as the corresponding asymptotic profile $u_M$, defined in \eqref{eq:selfsimburgers}. As we can observe, the solution given by \eqref{eq:scheme} is already quite close to $u_M$. However, a non-suitable viscous numerical flux like, for instance, the modified Lax-Friedrichs (e.g. \cite[Chapter 3]{Godlewski:1991}) can definitely modify the large-time behavior of the solution. In fact, in this case a viscosity proportional to $\dx^2/\dt$ is being added to the equation of the asymptotic profile (see \cite{Pozo:2015}), producing a more diffused wave. Nevertheless, the discretization of the non-linear term is not the only one with the ability to perturb the dynamics of the model. Let us emphasize that an inappropriate discretization of the non-local term also leads to an incorrect asymptotic profile. Note that in Figure \ref{fig1} we have the same scheme \eqref{eq:scheme} but taking $\coefzero=\coefone =1$, which produces a translated solution.

\begin{figure}[ht!]
	\includegraphics[width=0.85\linewidth]{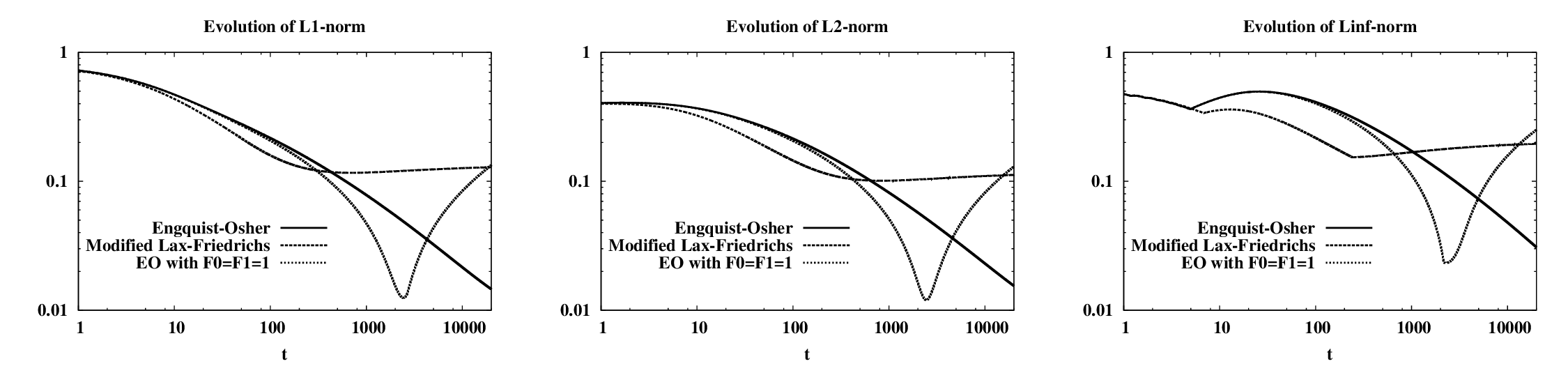}
	\caption{Evolution of the norms of the difference between the asymptotic profile and the solutions, multiplied by their corresponding rate $t^{\frac12(1-\frac1p)}$. From left to right, $L^1$, $L^2$ and $L^\infty$ norms. We compare \eqref{eq:scheme} (solid), modified LF numerical flux (dashed) and  EO with $\coefzero=\coefone=1$ (dotted).}
	\label{fig2}
\end{figure}

The convergence rates, given in \eqref{eq:convergencerates}, are shown in Figure \ref{fig2}. The graphic highlights the different performances mentioned above. In fact, the solution given by \eqref{eq:scheme} is the only one for which the norm is converging to zero with the corresponding rates.

\begin{figure}[ht!]
	\includegraphics[width=0.5\linewidth]{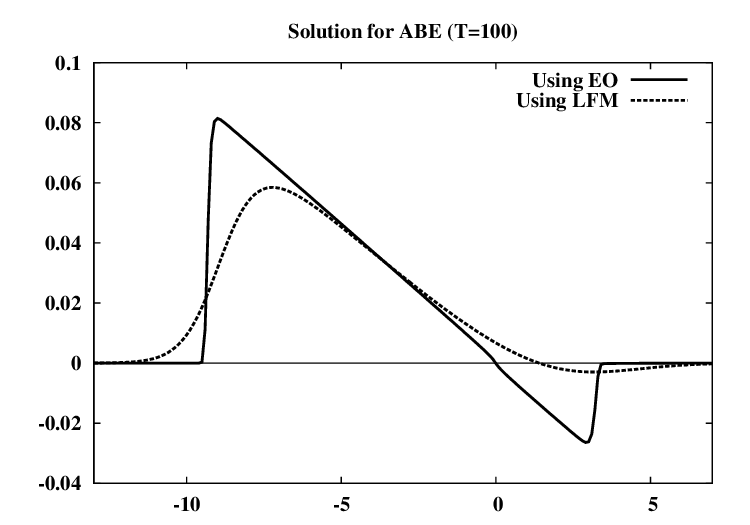}
	\caption{Solution of ABE with $\nu=10^{-4}$, $c=2\times10^{-4}$ and $\theta=1$ at $t=100$, using scheme \eqref{eq:scheme} discretized explicitly. We use Engquist-Osher (solid) and modified Lax-Friedrichs (dashed) numerical fluxes for the nonlinearity.}
	\label{fig3}
\end{figure}

To conclude, let us remark again the importance of taking a well-behaving numerical flux. In this paper we have proved that the asymptotic profile of \eqref{eq:model} is a diffusive wave. Therefore, any sign-changing initial data will lose its positive or negative part, depending on the sign of its mass. As in the case of the viscous Burgers equation \cite{Kim:2001}, simulations show that N-waves are intermediate states. Therefore, if the numerical viscosity is sufficiently large, the diffusion will become dominant much earlier than in the continuous model and destroy these profiles. For instance, let us consider the case $\nu=10^{-4}$ and $c=2\times10^{-4}$. In Figure \ref{fig3}, we can observe that at $t=100$ the N-wave shape is not preserved if the modified Lax-Friedrichs flux is used, while Engquist-Osher is able to keep the continuous dynamics. This numerical phenomenon was already observed in \cite{Pozo:2015} in the context of scalar conservation laws. It is interesting to see that the same pathology persist for viscous flows.

%%%%%%% APPENDIX %%%%%%%%%%%%
\appendix
\section{Auxiliary results} \label{appendix}
Here we prove some of the auxiliary results that we have used along the paper.

\begin{lmm}\label{lemma:int1}
For any piecewise constant function $w$ defined as in \eqref{eq:udelta} and $\dx>0$, the following holds:
\begin{equation*}
	\|w\|_p^{p(p+1)/(p-1)} \le 4 \|w\|_1^{2p/(p-1)}\|d^+_{\dx} |w|^{p/2}\|_2^2
\end{equation*}
for all $p\in(1,\infty)$.
\end{lmm}
\begin{proof}
First, let us define a piecewise linear function $v$ as follows:
\begin{equation*}
	v(x):=w_j\frac{x_{j+1}-x}{\dx}+w_{j+1}\frac{x-x_j}{\dx},\quad x\in[x_j,x_{j+1}].
\end{equation*}
On the one hand, we know that
\begin{equation}\label{eq:est3}
	\|v\|_\infty^2\le2\|v\|_2\|v_x\|_2.
\end{equation}
On the other hand, we have that:
\begin{equation*}
	\|v\|_{2}^{2}=\dx\sum_{j\in\Z}\int_{0}^{1}\left|w_j(1- x)+w_{j+1} x \right|^2dx \le \frac{1}{2} \dx\sum_{j\in\Z}\left(|w_j|^2+|w_{j+1}|^2\right) = \|w\|_2^2.
\end{equation*}
Moreover, it is easy to see that $\|v_x\|_2= \|d^+_{\dx} w\|_2$. Therefore, we can obtain a similar inequality as \eqref{eq:est3} for $w$:
\begin{equation}\label{eq:lemmapcest}
	\|w\|_\infty^2=\|v\|_\infty^2\le2\|v\|_2\|v_x\|_2\le 2\|w\|_2 \|d^+_{\dx} w\|_2.
\end{equation}
Applying this inequality to $|w|^{p/2}$, we deduce:
\begin{equation*}
	\| w \|_\infty^{2p}=\||w|^{p/2}\|_\infty^4 \le 4 \||w|^{p/2}\|_2^2 \|d^+_{\dx} |w|^{p/2}\|_2^2 = 4 \| w \|_p^p \|d^+_{\dx} |w|^{p/2}\|_2^2.
\end{equation*}
Thus, combining this with
\begin{equation*}
	\|w\|_p^{2p^2/(p-1)} \le \|w\|_\infty^{2p}\|w\|_1^{2p/(p-1)},
\end{equation*}
we conclude
\begin{equation*}
	\|w\|_p^{p(p+1)/(p-1)} \le 4 \|w\|_1^{2p/(p-1)}\|d^+_{\dx} |w|^{p/2}\|_2^2.
\end{equation*}
\end{proof}

\begin{proof}[Proof of Lemma \ref{lemma:conservation}]
For the first assertion, we simply integrate \eqref{eq:umu} over the whole space domain. We observe that all terms on the right hand side vanish, so 
\begin{equation*}
	 \frac{d}{dt}\int_\R \umu(t,x)dx = 0,\quad \forall t\ge0,
\end{equation*}
for all $\mu>0$ and, hence, the mass is conserved. Using the definition of $\umu$, we conclude
\begin{equation*}
	\int_\R \umu(t,x)dx = \int_\R \umu(0,x)dx = \int_\R \mu \ud^0(\mu x)dx = \int_\R \ud^0(x)dx.
\end{equation*} 

For the contractivity we prove that for any $u_0,v_0\in L^1(\R)$, their corresponding solutions $\umu$ and $\vmu$ satisfy
\begin{equation}\label{eq:l1contrac}
	\|\umu-\vmu\|_1 \le \|\umu_0-\vmu_0\|_1.
\end{equation}
For the sake of clarity, let us define $\wmu=\umu-\vmu$. Clearly, $\wmu$ verifies
\begin{align*}
		\wmu_t(t,x)&= \frac14 \Big(d^+_{\dx/\mu} \big(\umu(t,x)^2\big)+d^-_{\dx/\mu} \big(\umu(t,x)^2\big) - d^+_{\dx/\mu} \big(\vmu(t,x)^2\big) - d^-_{\dx/\mu} \big(\vmu(t,x)^2\big)\Big) \\[10pt]
		&\quad + \dx\, d^+_{\dx/\mu} R\big(\umu(t,x-\frac{\dx}{\mu}),\umu(t,x)\big) - \dx\, d^+_{\dx/\mu} R\big(\vmu(t,x-\frac{\dx}{\mu}),\vmu(t,x)\big) \\[5pt]
		&\quad + d^-_{\dx/\mu} \left(d^+_{\dx/\mu} \wmu(t,x) \right) \\[5pt]
		&\quad + \mu^2 \sum_{m=1}^{N} \omega_m \wmu(t,x-m\frac{\dx}{\mu}) -\mu^2 \coefzero \wmu(t,x) +\mu \coefone d^+_{\dx/\mu} \wmu(t,x).
\end{align*}
We multiply it by $\sign(\wmu)$ and integrate it on all $\R$. Using the definition of $R$ in \eqref{eq:numviscR} and reordering the terms we get
\begin{align}
	\frac{d}{dt} &\int_\R |\wmu(x)|dx \label{eq:contracdelta} \\
	&=\frac{1}{4} \int_\R d^+_{\dx/\mu}\big(\umu(x)^2+\umu(x-\frac{\dx}{\mu})^2+\umu(x)|\umu(x)|-\umu(x-\frac{\dx}{\mu})|\umu(x-\frac{\dx}{\mu})|\big) \sign(\wmu(x)) dx \notag \\
	&\quad - \frac{1}{4} \int_\R d^+_{\dx/\mu}\big(\vmu(x)^2+\vmu(x-\frac{\dx}{\mu})^2+\vmu(x)|\vmu(x)|-\vmu(x-\frac{\dx}{\mu})|\vmu(x-\frac{\dx}{\mu})|\big) \sign(\wmu(x)) dx \notag \\
	&\quad + \int_\R d^-_{\dx/\mu}(d^+_{\dx/\mu}(\wmu(x))) \sign(\wmu(x)) dx \notag \\
	&\quad + \mu^2 \sum_{m=1}^{N} \omega_m \int_\R \Big(\wmu(x-m\frac{\dx}{\mu}) -\wmu (x) \Big) \sign(\wmu(x)) dx +\mu \coefone \int_\R d^+_{\dx/\mu} \wmu(t,x) \sign(\wmu(x)) dx \notag \\
	&= I_1 + I_2 + I_3 + I_4 + I_5. \notag 
\end{align}
For $i=0,1$, let us denote $W^\pm_i=\{x\in\R:\pm \wmu(x-i\dx)>0\}$ and $W^0_i=\{x\in\R:\wmu(x-i\dx)=0\}$. Now we can split the domains of the integrals into several parts, according to the sign of $\wmu$. On the one hand, we have:
\begin{align*}
	I_1+I_2&= - \frac{1}{4} \int_\R \big(\umu(x)^2+\umu(x)|\umu(x)|\big) d^-_{\dx/\mu} (\sign(\wmu(x))) dx \\
	&\quad - \frac{1}{4} \int_\R \big(\umu(x-\frac{\dx}{\mu})^2-\umu(x-\frac{\dx}{\mu})|\umu(x-\frac{\dx}{\mu})|\big) d^-_{\dx/\mu} (\sign(\wmu(x))) dx \\
	&\quad + \frac{1}{4} \int_\R \big(\vmu(x)^2+\vmu(x)|\vmu(x)|\big) d^-_{\dx/\mu} (\sign(\wmu(x))) dx \\
	&\quad + \frac{1}{4} \int_\R \big(\vmu(x-\frac{\dx}{\mu})^2-\vmu(x-\frac{\dx}{\mu})|\vmu(x-\frac{\dx}{\mu})|\big) d^-_{\dx/\mu} (\sign(\wmu(x))) dx \\
	&= - \frac{\mu}{2\dx} \int_{ W_0^-\cap W_1^+} \big(\umu(x)^2+\umu(x)|\umu(x)|-\vmu(x)^2-\vmu(x)|\vmu(x)|\big) \sign(\wmu(x)) dx \\
	&\quad - \frac{\mu}{2\dx} \int_{ W_0^-\cap W_1^+} \big(\vmu(x-\frac{\dx}{\mu})^2-\vmu(x-\frac{\dx}{\mu})|\vmu(x-\frac{\dx}{\mu})| \\
	&\quad \qquad\qquad\qquad\qquad\qquad -\umu(x-\frac{\dx}{\mu})^2+\umu(x-\frac{\dx}{\mu})|\umu(x-\frac{\dx}{\mu})|\big) \sign(\wmu(x-\frac{\dx}{\mu})) dx \\
	&\quad - \frac{\mu}{2\dx} \int_{ W_0^+\cap W_1^-} \big(\umu(x)^2+\umu(x)|\umu(x)|-\vmu(x)^2-\vmu(x)|\vmu(x)|\big) \sign(\wmu(x)) dx \\
	&\quad - \frac{\mu}{2\dx} \int_{ W_0^+\cap W_1^-} \big(\vmu(x-\frac{\dx}{\mu})^2-\vmu(x-\frac{\dx}{\mu})|\vmu(x-\frac{\dx}{\mu})| \\
	&\quad \qquad\qquad\qquad\qquad\qquad -\umu(x-\frac{\dx}{\mu})^2+\umu(x-\frac{\dx}{\mu})|\umu(x-\frac{\dx}{\mu})|\big) \sign(\wmu(x-\frac{\dx}{\mu})) dx \\
	&\quad - \frac{\mu}{4\dx} \int_{ W_1^0} \big(\umu(x)^2+\umu(x)|\umu(x)|-\vmu(x)^2-\vmu(x)|\vmu(x)|\big) \sign(\wmu(x)) dx \\
	&\quad - \frac{\mu}{4\dx} \int_{ W_0^0} \big(\vmu(x-\frac{\dx}{\mu})^2-\vmu(x-\frac{\dx}{\mu})|\vmu(x-\frac{\dx}{\mu})|\\
	&\quad \qquad\qquad\qquad\qquad\qquad -\umu(x-\frac{\dx}{\mu})^2+\umu(x-\frac{\dx}{\mu})|\umu(x-\frac{\dx}{\mu})|\big) \sign(\wmu(x-\frac{\dx}{\mu})) dx.
\end{align*}
Using that
\begin{equation*}
	\big(b(b+|b|)-a(a+|a|)\big) \sign(b-a)\ge 0,\quad\forall a,b\in\R,
\end{equation*}
and that
\begin{equation*}
	\big(a(a-|a|)-b(b-|b|)\big) \sign(b-a)\ge 0,\quad\forall a,b\in\R,
\end{equation*}
we conclude that $I_1+I_2\le0$. On the other hand, since
\begin{align*}
	\int_\R \wmu(x-m\frac{\dx}{\mu}) \sign(\wmu(x)) dx \le \int_\R |\wmu (x)| dx, \quad\forall m\in\Z,
\end{align*}
it is immediate that 
\begin{align*}
	I_3 &= \frac{\mu^2}{\dx^2} \int_\R \Big( \wmu(x-\frac{\dx}{\mu}) + \wmu(x+\frac{\dx}{\mu})  - 2\wmu(x) \Big)\sign(\wmu(x)) dx \le0.
\end{align*}
Moreover, for the same reason, we deduce that $I_4\le$ and $I_5\le 0$. Therefore, from \eqref{eq:contracdelta} we get that
\begin{align}
	\frac{d}{dt} &\int_\R |\wmu(x)|dx \le 0,
\end{align}
This guarantees the contractive property \eqref{eq:l1contrac}.
\end{proof}

\begin{proof}[Proof of Lemma \ref{lemma:hs}]
Let us consider the Fourier transform of $w$ as
\begin{equation*}
	\widehat{w}(\xi)=\int_\R e^{-ix\xi}w(x)dx,\quad\xi\in\R
\end{equation*}
and the discrete Fourier transform of the sequence $\{w_j\}_{j\in\Z}$ as
\begin{equation*}
	\squarehat{w}(\xi)=\dx\sum_{j\in\Z} w_j e^{-ij\dx\xi},\quad\xi\in[-\frac{\pi}{\dx},\frac{\pi}{\dx}].
\end{equation*}
It is also clear that for a piecewise constant function $w$ defined as in \eqref{eq:udelta}
\begin{equation*}
	\widehat{w}(\xi)=\frac{2\sin(\frac{\xi\dx}{2})}{\xi\dx} \squarehat{w}(\xi) \quad \mbox{ and } \quad \widehat{d^+_{\dx}w}(\xi)=\frac{e^{i\xi\dx}-1}{\dx}\widehat{w}(\xi).
\end{equation*}
Now, we know that
\begin{equation}\label{eq:normds}
	\||D|^sw\|_2^2=\int_\R |\xi|^{2s}|\widehat{w}(\xi)|^2d\xi =\int_{-\pi/\dx}^{\pi/\dx} |\xi|^{2s}|\widehat{w}(\xi)|^2d\xi + \sum_{j\ne0} \int_{(2j-1)\pi/\dx}^{(2j+1)\pi/\dx} |\xi|^{2s}|\widehat{w}(\xi)|^2d\xi.
\end{equation}
For each $j\ne0$, we have
\begin{align*}
	\int_{(2j-1)\pi/\dx}^{(2j+1)\pi/\dx} |\xi|^{2s}|\widehat{w}(\xi)|^2d\xi & = \int_{(2j-1)\pi/\dx}^{(2j+1)\pi/\dx} |\xi|^{2s}|\squarehat{w}(\xi)|^2\left|\frac{2\sin(\frac{\xi\dx}{2})}{\xi\dx}\right|^2d\xi \\
	& = \int_{-\pi/\dx}^{\pi/\dx} \left|\xi+2j\frac{\pi}{\dx}\right|^{2s}|\squarehat{w}(\xi)|^2 \left|\frac{2\sin(\frac{\xi\dx}{2}+j\pi)}{(\xi+2j\frac{\pi}{\dx})\dx}\right|^2 d\xi \\
	& = \int_{-\pi/\dx}^{\pi/\dx} \left|\frac{2}{\dx}\sin(\frac{\xi\dx}{2})\right|^{2s} |\squarehat{w}(\xi)|^2 \left|\frac{2}{\dx}\sin(\frac{\xi\dx}{2})\right|^{2-2s} \left|\xi+2j\frac{\pi}{\dx}\right|^{2s-2} d\xi \\
	& \le \int_{-\pi/\dx}^{\pi/\dx} |\xi|^{2s} |\squarehat{w}(\xi)|^2 \frac{|\xi|^{2-2s}}{\left|\xi+2j\frac{\pi}{\dx}\right|^{2-2s}} d\xi \\
	& \le \frac{1}{\left|2|j|-1\right|^{2-2s}} \int_{-\pi/\dx}^{\pi/\dx} |\xi|^{2s} |\squarehat{w}(\xi)|^2 d\xi.
\end{align*}
Therefore, replacing this in \eqref{eq:normds} and using that $0<s<\frac12$, we get
\begin{align*}
	\||D|^sw\|_2^2&\le\int_{-\pi/\dx}^{\pi/\dx} |\xi|^{2s}|\widehat{w}(\xi)|^2d\xi + \sum_{j\ne0} \frac{1}{\left|2|j|-1\right|^{2-2s}} \int_{-\pi/\dx}^{\pi/\dx} |\xi|^{2s} |\squarehat{w}(\xi)|^2 d\xi \\
	&= \int_{-\pi/\dx}^{\pi/\dx} |\xi|^{2s} \left|\frac{2\sin(\frac{\xi\dx}{2})}{\xi\dx}\right|^2 |\squarehat{w}(\xi)|^2d\xi + \sum_{j\ne0} \frac{1}{\left|2|j|-1\right|^{2-2s}} \int_{-\pi/\dx}^{\pi/\dx} |\xi|^{2s} |\squarehat{w}(\xi)|^2 d\xi \\
	&\lesssim \int_{-\pi/\dx}^{\pi/\dx} |\xi|^{2s} |\squarehat{w}(\xi)|^2d\xi.
\end{align*}
On the other hand, using analogous arguments, we also have
\begin{align*}
	\|d^+_{\dx}w\|_2^2 &= \int_\R \left|\frac{e^{i\xi\dx}-1}{\dx}\right|^2|\widehat{w}(\xi)|^2d\xi = \int_\R \left|\frac{e^{i\xi\dx}-1}{\dx}\right|^2 \left|\frac{2\sin(\frac{\xi\dx}{2})}{\xi\dx}\right|^2 |\squarehat{w}(\xi)|^2d\xi \\
	&\gtrsim \int_{-\pi/\dx}^{\pi/\dx} \left|\frac{e^{i\xi\dx}-1}{\dx}\right|^2 |\squarehat{w}(\xi)|^2d\xi \gtrsim \int_{-\pi/\dx}^{\pi/\dx} \left|\xi\right|^2 |\squarehat{w}(\xi)|^2d\xi
\end{align*}
Finally, we conclude
\begin{align*}
	\|w\|_{H^s(\R)}^2&=\int_\R(1+|\xi|^{2s})|\widehat w(\xi)|^2d\xi =\int_\R |\widehat w(\xi)|^2d\xi+\int_\R |\xi|^{2s}|\widehat w(\xi)|^2d\xi \\
	& \lesssim \int_{-\pi/\dx}^{\pi/\dx} |\squarehat w(\xi)|^2d\xi + \int_{-\pi/\dx}^{\pi/\dx} |\xi|^{2s}|\squarehat{w}(\xi)|^2 d\xi \lesssim \int_{-\pi/\dx}^{\pi/\dx} (1+|\xi|^{2s})|\squarehat{w}(\xi)|^2 d\xi \\
	& \lesssim \int_{-\pi/\dx}^{\pi/\dx} (1+|\xi|^{2})|\squarehat{w}(\xi)|^2 d\xi \lesssim (\|w\|_2^2+\|d^+_{\dx} w\|_2^2).
\end{align*}
\end{proof}

\begin{lmm}\label{lemma:serie}
Given any $a\in(0,1)$ and $b\in\C$ with $|b|=1$, the following inequality holds:
\begin{equation*}
	\left|\sum_{m=1}^N a^{m} (b^m-1) + \left(\sum_{m=1}^N m a^{m}\right)\left(\frac{1}{b}-1\right)\right| \le \left|b-1\right|^2 \frac{a}{(1-a)^3}.
\end{equation*}
\end{lmm}
\begin{proof}
Using that $|b|=1$, we have:
\begin{align*}
	\left|\sum_{m=1}^N a^{m} (b^m-1) \right.&\left.+\left(\sum_{m=1}^N m a^{m}\right)\left(\frac{1}{b}-1\right)\right| =\left|\sum_{m=1}^N a^{m} b(b^m-1) - \left(\sum_{m=1}^N m a^{m}\right)\left(b-1\right)\right| \\
	&=\left|b-1\right|\left|\sum_{m=1}^N a^{m} \sum_{k=0}^{m-1}(b^{k+1}-1)\right| =\left|b-1\right|^2\left|\sum_{m=1}^N a^{m} \sum_{k=0}^{m-1}\sum_{j=0}^{k} b^j \right| \\
	& \le \frac{1}{2}\left|b-1\right|^2 \sum_{m=1}^\infty m(m+1)a^{m}=\left|b-1\right|^2 \frac{a}{(1-a)^3}.
\end{align*}
\end{proof}

\section*{Acknowledgements}
The authors would like to thank Enrique Zuazua (DeustoTech and UAM) for the stimulating discussions that led to the questions addressed in this paper. This work was partially supported by the Advanced Grants NUMERIWAVES/FP7-246775 of the European Research Council Executive Agency, PI2010-04 program of the Basque Government, the Grant MTM2014-52347 of Spanish Ministry of Economy and Competitiveness and also by the Basque Government through the BERC 2014-2017 program and by Spanish Ministry of Economy and Competitiveness MINECO: BCAM Severo Ochoa excellence accreditation SEV-2013-0323. L.~I.~Ignat was also partially supported by Grant  PN-II-ID-PCE-2011-3-0075  of the Romanian National Authority for Scientific Research, CNCS-UEFISCDI and by FA9550-15-1-0027 of AFOSR. A.~Pozo was granted by the Basque Government, reference PRE\_2013\_2\_150.

\bibliography{library}
\bibliographystyle{amsplain}

\end{document}